\input graphicx

\def \WB {\widetilde B}
\def\tttt #1{{\textstyle{#1} }}
\def \ul #1  {{\underline{#1} }}
\def \ol #1  {{\overline{#1} }}

\def \bu {{\hskip -.1in}}
\def \DD {\Delta}

\def\con {\subseteq}

\def \CF {{\cal F}}

\def\la{{\lambda}}

\def \CD {{\cal D}}

\overfullrule=0pt
\baselineskip 14pt
\settabs 8 \columns
\+ \hfill&\hfill \hfill & \hfill \hfill   & \hfill \hfill  & \hfill \hfill   &
 & \hfill \hfill\cr
\parindent=.5truein
\hfuzz=3.44182pt
\hsize 6.5truein

\font\ita=cmssi10 
 \font\small=cmr6
\font\title=cmbx10 scaled\magstep2
\font\normal=cmr10 
\font\small=cmr6

\font\bol=cmbx12
\def\del {\partial}
\def \con {\subseteq}

\def\sig{\sigma}

\def \-> {\rightarrow}
\def\LL{\big\langle}

\def\RR {\big\rangle}

\def\DD {\Delta}

\def\om {\omega}
\def\la {\lambda}

\def \RA {\rightarrow}

\def \sas {\vskip .06truein}
\def\sa{{\vskip .125truein}}

\def\sap{{\vskip .25truein}}

\def \eee {\epsilon}

\def\aaa {\alpha}

\def\aa {\alpha}
\def\bb {\beta}

\def\gg {\gamma}
\def\con {\subseteq}
\def \ses {\enskip = \enskip}
\def \sps {\enskip + \enskip}

\def \sms {\enskip -\enskip}

\def \scs {\ssp , \ssp}
\def \ess {\enskip}
\def \ssp {\hskip .25em}
\def \bigsp {\hskip .5truein}
\def \part {\vdash}

\def \DD {\Delta}

\normal

\vsize=9truein
\sap
\def\today{\ifcase\month\or
January\or February\or March\or April\or may\or June\or
July\or August\or September\or October\or November\or
December\fi
\space\number\day, \number\year}

\headline={   
\small \today \hfill \hfill$\ess\ess\ess$$\ess\ess\ess\ess\ess\ess$ \folio }
 \footline={\hfil}

\def \RA {{ \rightarrow }}

\def \BN {{\bf N}}

\def \BQ {{\bf Q}}

\def \II {{\rm I}}
\def \om {\omega}

\def \TH {{\widetilde H}}

\def \BQ {{\bf Q}}
\def \om {\omega}

\def \TH {{\widetilde H}}

\font\small=cmr6
\def \scs {\ssp , \ssp}
\def \ess {\enskip}
\def \ssp {\hskip .25em}
\def \bigsp {\hskip .5truein}
\def \part {\vdash}
\font\title=cmbx10 scaled\magstep2
\font\normal=cmr10 
\def\today{\ifcase\month\or
January\or February\or March\or April\or may\or June\or
July\or August\or September\or October\or November\or
December\fi
\space\number\day, \number\year}
\ 
\vskip -.5in
\centerline{\title  e-Positivity Results and  Conjectures} 
\sa

\centerline{\bf A. M. Garsia, J. Haglund, D. Qiu   and M. Romero  \footnote{$(\dagger)$}{Research supported by NSF, the first and third  author by grant DMS-1700233, the second by grant DMS-1600670  the last  by grant DMS-1362160.}}
\sas

\noindent{\bol Abstract} 

In a 2016 ArXiv posting F. Bergeron listed a variety of symmetric functions $G[X;q]$ with the property that $G[X;1+q]$ is $e$-positive.  A large subvariety of his examples could be  explained by the conjecture that the  Dyck path LLT  polynomials exhibit the same phenomenon. In this paper we list the results of computer explorations which suggest that other examples exhibit the same phenomenon. We prove two of the resulting conjectures and propose algorithms that would prove several
of our conjectures. In writing this paper we have learned that similar findings have been independently discovered by Per Alexandersson (see [1]).
\sas

\noindent{\bol Introduction} 
\sas

We say that the symmetric function $G[X;q]$
exhibits the $e$-positivity phenomenon if and only if the symmetric function
 $G[X;1+q]$ is $e$-positive.
 This only means that, in the $e$-basis expansion
$$
G[X;1+q]\ses \sum_{\la}
a_\la(q) \, e_\la[X],
\eqno \II.1
$$
the coefficients $a_\la(q)$ are polynomials in $q$ with positive integer coefficients.
The following are four examples:
$$
LLT(4,3,[0,1,2])[X;1+q]\ses(q^3+2 q^2) e_3+(q^2+3 q) e_2 e_1+
e_1^3,
\eqno \II.2
$$
$$
LLT(7,4,[0,1,2,2])[X;1+q]\ses
(q^3+2 q^2) e_4+(q^2+2 q) e_3 e_1 +q e_2^2
+ e_2e_1^2,
\eqno \II.3
$$
$$
B_{[3,1,1]}[X;1+q]\ses
(q^3+2 q^2) e_5+(q^2+2 q) e_4 e_1+
q e_3e_2 +e_3e_1^2,
\eqno \II.4
$$
$$
Unicell_{[1,4,3,2]}[X;q]\ses(s_{[4]}+(2q+1) s_{[3, 1]}+2 q s_{[2,2] }+(q^2+2 q) s_{[2, 1, 1]}+q^2 s_{[1, 1, 1, 1]},
\eqno \II.5
$$
$$
Unicell_{[1,4,3,2]}[X;1+q]\ses q^2 e_3e_1 +2 e_2e_1^2  q+e_1^4.
\eqno \II.6
$$
The first is the LLT polynomial of the path that alternates North steps and East steps,
the second is the LLT of a rational Dyck path in the $7\times 4$ lattice rectangle. The third is a balanced Dyck  path LLT that hits the diagonal according to the partition 
$[3,1,1]$. In I.5 we have the unicellular LLT  whose successive cells  are in diagonals $1,4,3,2$. In I.6 we see that even in the latter case the $e$-positivity phenomenon takes place.

The  experimental evidence of widest impact we have noticed  so far is that the LLT polynomials generated by Dyck paths whether classical or rational do exhibit the $e$-positivity phenomenon.

Since a Dyck path that alternates North and East steps is also Balanced we tested some 
Balanced paths and sure enough, the  $e$-positivity phenomenon seems to occur there as well.   Noticing that N-E alternating 
Dyck paths are also unicellular, we tested several cases  of these  LLT's and discovered that the  $e$-positivity phenomenon seems to occur there too. Our experimental data lead us to conjecture that the following  
families of symmetric functions exhibit the $e$-positivity phenomenon.
\sas

{$\bf (1)$} {\ita The Modified Macdonald polynomials at $t=1$. That is $\TH_\mu[X;q,1]$, for any partition $\mu$.}  

{$\bf (2)$} {\ita The  polynomials $B_p[X;q]$ for all compositions $p$ (see section 2).}

{$\bf (3)$} {\ita All unicellular LLT polynomials.} 

{$\bf (4)$} {\ita All column LLT polynomials, see 
the last section of this paper.}

{$\bf (5)$}
{\ita The polynomials $\nabla C_p\, 1$ for all compositions $p$.}

 {$\bf (6)$} {\ita All the  polynomials  $Q_{m,n}(-1)^n$, appearing in the rational Shuffle  Conjecture.}
 
 {$\bf (7)$} {\ita All the  polynomials
 $\DD_{e_k}e_n$ appearing in the Delta conjecture} 

{$\bf (8)$} {\ita All the polynomials $\nabla (-1)^{n-1}p_n$.}

\vfill
\supereject

Our main results here are proofs of conjectures 
$\bf (1)$ and $\bf (2)$ and the construction of algorithms  that would yield the asserted $e$-basis expansions for classical Dyck path, unicellular and column LLT's. 

Our presentation is divided into three sections. In the first section we prove  conjecture $\bf (1)$. In the second  section we 
prove conjecture $\bf (2)$. In the third and final section  we comment on some consequences of our results and state our conjectured $e$-expansion formula for column LLT polynomials.

In fact, it will be seen that, given a classical Dyck path $D$ and its zeta image $\zeta(D)$, by means of the areaprime way of computing $LLT's$, we can deal with classical Dyck path LLT's, column LLT's and unicellular LLT's at the same time by simply marking, partially marking and not marking the removable corners of the English partition above 
$\zeta(D)$. This view point   makes evident   that the number of distinct  polynomials of unicellular and column LLT's are respectively not larger than the Catalan numbers and lower Schr\"oder numbers.
We also show how to use the Carlsson-Mellit  algorithm for constructing these  LLT's to confirm  our conjectures for larger scale examples than what is achieved by purely combinatorial means.
\sa

\noindent{\bol 1. Proof of conjecture (1)}

For notation and plethystic notation we refer to [4] and [12] where the reader can also consult a   Modified Macdonald polynomials ``tool kit''.
Our point of departure is the following specialization at $t=1$ of the Modified Macdonald polynomial. 
\sas

\noindent{\bol Proposition 1.1}

{\ita For any partition $\mu$ we have
$$ 
\TH_\mu[X;q,1]\ses \prod_{i=1}^{\l(\mu)}
(q;q)_{\mu_i} 
h_{\mu_i}\big[\tttt{X\over 1-q}\big],
\eqno 1.1
$$
where for any integer $k\ge 0$ we have}
 $
\,\,(q;q)_k=(1-q)(1-q^2)\cdots (1-q^k).
 $
 
\noindent{\bol Proof}

In [16]  Chapter 8. Integral Forms  (see (8.4) Remark  (iii)) Macdonald proves that
$$
J_{\mu }[X;1,t]\ses (t;t)_{\mu' }e_{\mu' }(X).
\eqno 1.2
$$
Our definition of the Modified Macdonald
polynomial indexed by $\mu$ is
$$
\TH_\mu[X,q,t]\ses t^{n(\mu)}J_\mu\big[\tttt{X\over 1-1/t};q,1/t\big].
\eqno 1.3
$$
Thus setting $q=1$ and using 1.2 gives
$$
\TH_\mu[X;1,t]\ses  t^{n(\mu)} (1/t\,;1/t)_{\mu' }e_{\mu' }\big[\tttt{X\over 1-1/t}\big]. \eqno 1.4
$$
Now for $\mu\part m$,
$$
\eqalign{
t^{n(\mu)} (1/t\,;1/t)_{\mu' }
&\ses
t^{n(\mu)}  
\prod_{i=1}^{\l(\mu')}\prod_{j=1}^{\mu_i'}(1-1/t^{j})    
\ses
t^{n(\mu)}
\prod_{i=1}^{\l(\mu')}
t^{-\mu_i'}
  t^{-{\mu_i'\choose 2}}  
\prod_{i=1}^{\l(\mu')}(-1)^{\mu_i'}
\prod_{j=1}^{\mu_i'}{ (1-t^j)     }
\cr
&\ses
(-t)^{-m}
\prod_{i=1}^{\l(\mu')}
\prod_{j=1}^{\mu_i'}{ (1-t^j)}
\ses
(-t)^{-m}
(t,t)_{\mu'},
\cr} 
\eqno 1.5
$$
and
$$
\eqalign{
e_{\mu' }\big[\tttt{X\over 1-1/t}\big]
\ses
t^m e_{\mu' }\big[-\tttt{X\over 1-t }\big]
\ses 
(-t)^m h_{\mu' }\big[ \tttt{X\over 1-t }\big].
\cr
}
\eqno 1.6
$$
Combining 1.6 with 1.5 and 1.4 gives
$$
\TH_\mu[X;1,t]\ses (t,t)_{\mu'}h_{\mu' }\big[ \tttt{X\over 1-t }\big].
$$
Using the identity
 $  \,
\TH_{\mu'}[X;q,1]=\TH_{\mu }[X;1,q] \, , 
  $
we finally derive that
$$
\TH_{\mu'}[X;q,1]\ses (q,q)_{\mu'}h_{\mu' }\big[ \tttt{X\over 1-q }\big].
\,\,$$
But this is just another way of writing 1.1.

\supereject

Since the dual of the $e$-basis with respect to the Hall scalar product is the forgotten basis, for any integer $m\ge 1$ we obtain the 
$e$-basis expansion
$$
h_m\big[\tttt{X\over 1-q} \big]\ses\sum_{\mu\part m}e_\mu[X] f_\mu\big[\tttt{1\over 1-q} \big].
\eqno 1.7
$$
This given, in view of 1.1, to show that
$$
\TH_\mu[X;1+q,1]\ses  \prod_{i=1}^{\l(\mu)}
(q;q)_{\mu_i} 
h_{\mu_i}\big[\tttt{X\over 1-q}\big]
\Big|_{q=1+q}
\eqno 1.8
$$
is  $e$-positive it is sufficient to prove
the $e$-positivity of the polynomial
$
(q;q)_m h_m\big[\tttt{X\over 1-q}\big]
\Big|_{q=1+q}
$
for every $m\ge 1$.
But that will be true if and only if 
 we have
$$
(q;q)_m  f_\mu\big[\tttt{1\over 1-q} \big]\Big|_{q=1+q}\in \BN[q]\ess\ess\ess \hbox{(for all $\mu\part m$)}.
\eqno 1.9
$$
Remarkably computer data revealed that this fact is due to the general validity of 
the following identity
\sas

\noindent{\bol Proposition 1.1}

{\ita For any partition $\mu\part m$ we have
$$
f_\mu\big[\tttt{1\over 1-q}\big](q;q)_m\ses
\Pi_\mu(q)\, (q-1)^{m-\l(\mu)}
$$
with $\Pi_\mu(q)\in \BN[q]$}.
\sa

Thus our final goal in this section will be the proof of this result. It develops that 
to do  this we need  auxiliary identities some of which are well known.
For sake of completeness, we will give complete proofs  of all the needed identities.
 
\vskip -.138in
\hfill
${\includegraphics[width=1.7in]{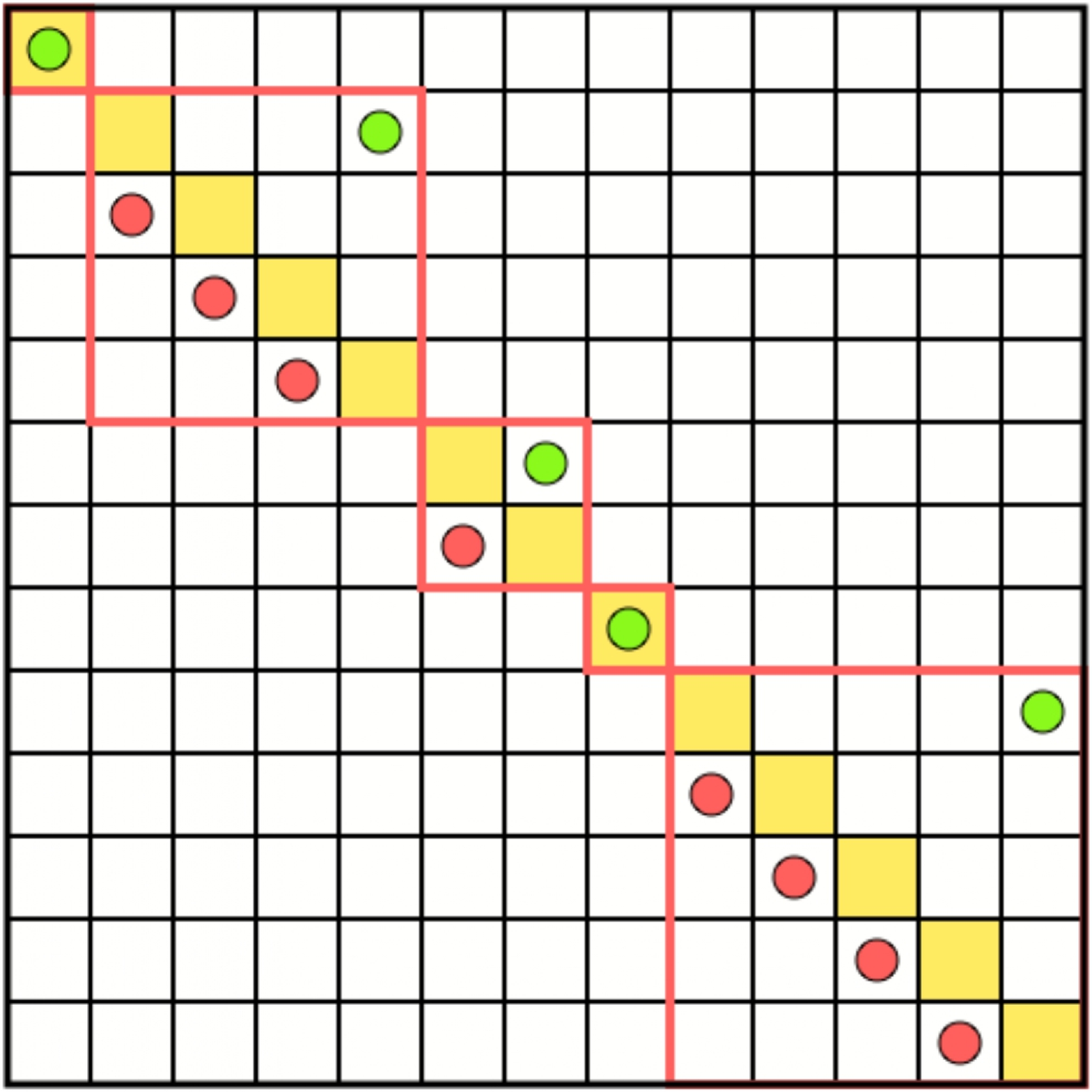}}  $

\hsize 4.6in
\vskip -1.55in
 We will start by dealing  with the factor $f_\mu\big[\tttt{1\over 1-q} \big]$. To this end recall that the Jacobi Trudi identity gives $h_m[X]=\det\| e_{j-i+1}\|_{i=1}^m$. Since this matrix has $1's$ in the subdiagonal and nothing but zeros below them, the only non vanishing determinantal terms are as indicated in  the adjacent figure.
We see there that  each subset of the sub-diagonal is broken up  into a union of strings of  adjacent 
elements. Each string determines a cycle  of the corresponding non vanishing term. The cycles are of the form
 $
(i,i+1,i+2,\ldots ,j)
 $
and  contribute  to the determinantal  
term  the factor
 $
(-1)^{j-i}\, e_{j-i+1}.
 $

\hsize 6.5in
\noindent

This example produces the 
term $(-1)^{13-5}e_1e_4e_2e_1e_5$. 
In particular, it follows that the coefficient of the $e$-basis element $e_5 e_4 e_2 e_1^2$ is equal to the number of distinct rearrangements
of the cycles yielding this product.
Thus the  general result may be written in the form
$$
h_m[X]=\sum_{\mu\part m}(-1)^{m-\l(\mu)}
|DR(\mu)| e_\mu[X],
$$
where $DR(\mu)$ is the set of distinct 
rearrangements $\aa_1,\aa_2,\ldots, \aa_{\l(\mu)}$ of $\mu_1,\mu_2,\ldots, \mu_{\l(\mu)}$. On the other hand we have
$$
h_m[X\cdot 1]\ses \sum_{\mu\part m}f_\mu{[1]} e_\mu[X].
$$
Therefore, we have for any monomial $\gg$,
$$
f_\mu[\gg]\ses \gg^m (-1)^{m-\l(\mu)}
|R(\mu)|\ses  (-1)^{m-\l(\mu)}
\sum_{\rho\in R(\mu)}\gg^{|\rho|}.
\eqno 1.10
$$
\supereject

To use this identity, we need  the following
\sas

\noindent{\bol Proposition 1.2 }
$$
f_\mu[X+Y]\ses \sum_{\aa\cup \bb=\mu}
f_\aa[X]f_\bb[Y],
\eqno 1.11
$$
{\ita where $\aa$ as well as $\bb$ are allowed to be empty partitions} 

\noindent{\bol Proof}

We have 
$$
\eqalign{
 e_n\big[X(Y+Z) \big]=\sum_{\mu\part n}
h_\mu[X]f_\mu[Y+Z] 
&=\sum_{k=0}^n
e_{n-k}\big[X Y\big]e_{ k}\big[X Z  \big]
=\sum_{k=0}^n
\sum_{\aa\part n-k}
\sum_{\bb\part k}
h_\aa[X]f_\aa[Y ]
h_\bb[X ]f_{ \bb}\big[ Z  \big]
\cr}.
$$
Equating coefficients of $h_\mu[X]$ gives 1.11.
\sas

For a sequence of partitions 
$\nu=(\nu^1,\nu^2,...)$ whose parts rearrange to $\mu$ we will write $\nu\in PR(\mu)$. Analogously, if 
$ p=( p^1,p^2,...)$ is a sequence of compositions whose parts rearrange to $\mu$
we will write $p\in CR(\mu)$.  In both cases we must allow the parts to be empty.
In particular, 1.11 may be rewritten in the form
$$
f_\mu[X+Y]\ses \sum_{(\aa,\bb)\in PR(\mu)}
f_\aa[X]f_\bb[Y].
\eqno 1.12
$$
Iterating this relation we obtain, for arbitrary $n$
$$
f_\mu[x_1+x_2+\cdots+x_{n} ]=\sum_{(\nu^1,\nu^2,\ldots,\,\nu^{n})\in PR(\mu)}\prod_{i= 1}^{n}f_{\nu^i}[x_i].
$$
Using 1.10 this  may be rewritten as
$$
\eqalign{
f_\mu[x_1+x_2+\cdots+ x_ n]
&\ses 
\sum_{(\nu^1,\nu^2,\ldots,\,\nu^n)\in PR(\mu)}\prod_{i= 1}^{n}
(-1)^{|\nu^i|-\l(\nu^i )}\sum_{\rho^i\in R(\nu^i)}x_i^{|\rho^i|}
\cr
&\ses
(-1)^{|\mu|-\l(\mu )}
\bu\bu
\sum_{(\nu^1,\nu^2,\ldots,\,\nu^{n})\in PR(\mu)}
\,\,\prod_{i= 1}^{n}
\sum_{\rho^i\in R(\nu^i)}x_i^{|\rho^i|}
\cr
&\ses
(-1)^{|\mu|-\l(\mu )}
\bu\bu
\sum_{p=(p_1,p_2,\ldots,p_{n})\in CR(\mu)}
\,\,  x_1^{|p_1|}x_2^{|p_2|}\cdots x_{n}^{|p_{n}|}.
\cr
}
$$
Now let $n\RA \infty$ to get
$$
f_\mu[x_1+x_2+x_3+\cdots]\ses
(-1)^{|\mu|-\l(\mu )}
\bu\bu
\bu\bu
\sum_{p=(p_1,p_2,p_3,\ldots)\in CR(\mu)}
\,\,  x_1^{|p_1|}x_2^{|p_2|}x_3^{|p_3|}\cdots .
$$
To  compute $f_\mu\big[\tttt{1\over 1-q}\big]$
we need only make the replacement $x_i\RA q^{i-1}$ obtaining
$$
f_\mu\big[\tttt{1\over 1-q}\big]\ses
(-1)^{|\mu|-\l(\mu )}
\bu\bu
\bu\bu
\sum_{p=(p_1,p_2,\ldots, p_{i},\ldots)\in CR(\mu)}
\,\,  (q^0)^{|p_1|}(q^1)^{|p_2|}\cdots (q^{i-1})^{|p_i|}\cdots .
\eqno 1.13 
$$
Now in the case $\mu=(4,3,3,2,1,1)$ one of the possible summands is
$$
p=\big( (1,3) \scs(2,4,1)\scs\phi\scs(3) \big)\in
CR((4,3,3,2,1,1)).
$$
The corresponding term in the sum is the monomial
$$
(q^0)^{|(1,3)|}(q^1)^{|(2,4,1)|}(q^2)^{|\phi|}(q^3)^{|(3)|}\, =\, (q^0)^1(q^0)^3\cdot (q^1)^2(q^1)^4(q^1)^1
\cdot (q^3)^3
$$

We clearly will obtain this case with the specialization $a=(1,3,2,4,1,3)$ and
$$
i_1=i_2=0\scs i_3=i_4=i_5=1\scs \ess\hbox{and}\ess\ess i_6=3
$$

\supereject
x in the sum
$$
f_{\mu}\big[\tttt{1\over 1-q}\big]\ses (-1)^{|\mu|-\l(\mu )}
\bu\bu\bu
\sum_{a= (a_1,a_2,\cdots,a_{\l(\mu)})\in DR(\mu)}\,\,\,
\sum_{0\le i_1\le i_2\le \cdots \le i_{\l(\mu)}}
\,\,\,(q^{a_1})^{i_1 }(q^{a_2})^{i_2 }
\cdots (q^{a_{\l(\mu)}})^{i_{\l(\mu)} }
\eqno 1.14
$$
for $\mu=(4,3,3,2,1,1)$. Here, again   $DR(\mu)$  denotes the collection of all distinct rearrangements of the components of $\mu$. Which for  our example it is a collection of ${6!\over 1! 2! 1! 2!}=180$
elements. In fact, a moment's reflection should reveal that, in full generality, the construction we gave for our particular example can be extended to obtain a bijection between the terms in the right hand sides of 1.13 and 1.14. This 
 assures the equality of the resulting rational functions of $q$.  

Now, by a very simple trick, we can obtain an explicit formula for the rational function 
$$
F_{(a_1,a_2,\ldots ,a_{\l(mu)})}(q)
\ses
\sum_{0\le i_1\le i_2\le \cdots \le i_{\l(\mu)}}
\,\,\,(q^{a_1})^{i_1 }(q^{a_2})^{i_2 }
\cdots (q^{a_{\l(\mu)}})^{i_{\l(\mu)} }.
\eqno 1.15
$$
The standard step is to simply  make the change of variables
$$
i_1=r_1,\ess i_2=r_1+r_2,\ess \ldots ,\ess
i_{\l(\mu)}=r_1+r_2 + \cdots +r_{\l(\mu)} 
$$ 
and rewrite 1.15 in the form
$$
\eqalign{
F_{(a_1,a_2,\ldots ,a_{\l(mu)})}(q)
&\ses
\sum_{r_1\ge 0}\sum_{r_2\ge 0}\cdots \sum_{r_{\l(\mu)}\ge 0}
\,\,\, (q^{a_1+a_2+\cdots +a_{\l(\mu)} })^{r_1 }
(q^{a_2 +\cdots +a_{\l(\mu)}} )^{r_2}
\cdots (q^{a_{\l(\mu)}})^{ r_{\l(\mu)} }
\cr
&\ses
{1\over 1- q^{a_1+a_2+\cdots +a_{\l(\mu)} }}
{1\over 1- q^{a_2+\cdots +a_{\l(\mu)} }}\cdots
{1\over 1- q^{a_{\l(\mu)} }} .
\cr
}
\eqno 1.16
$$
Thus from 1.14 it follows that
$$
(q;q)_m \, f_{\mu}\big[\tttt{1\over 1-q}\big]= (-1)^{|\mu|-\l(\mu )}
\bu\bu\bu\bu 
\sum_{a= (a_1,a_2,\cdots,a_{\l(\mu)})\in DR(\mu)}
{(1-q)(1-q^2)(1-q^3)\cdots (1-q^m) \over \big(1- q^{a_1+a_2+\cdots +a_{\l(\mu)} }\big)
\big(1- q^{a_2+\cdots +a_{\l(\mu)} }\big)
 \cdots  \big(1- q^{a_{\l(\mu)}}\big)
 }.
\eqno 1.17
$$
This given we are now ready to give our
\sas
\noindent{\bol Proof of Proposition 1.1}
\sas
Since the components of $a=(a_1,a_2,\cdots ,a_{\l(\mu)})$ are only a rearrangement of the components of $\mu\part m$, the integers
$$
a_1+a_2+a_3+\cdots +a_{\l(\mu)}, \ess
 a_2+a_3+\cdots +a_{\l(\mu)}  
 , \ess\cdots, \ess
 a_{\l(\mu)}
$$
are distinct and therefore form a subset of $\{1,2,3,\ldots, m\}$. Let us then set
$$
S(a)\ses \{1,2,3,\ldots, m\}\, -\, \{ a_1+a_2+\cdots +a_i\, :\, i=1,2,\ldots , \l(\mu)\}.
\eqno 1.18
$$
Thus we can rewrite  1.17 as
$$
\eqalign{
(q;q)_m \, f_{\mu}\big[\tttt{1\over 1-q}\big]&\ses (-1)^{|\mu|-\l(\mu )}(1-q)^{m-\l(\mu)}
\bu\bu\bu\bu 
\sum_{a= (a_1,a_2,\cdots,a_{\l(\mu)})\in DR(\mu)}\prod_{i\in S(a)}[i]_q 
\cr
&\ses
(q-1)^{|\mu|-\l(\mu )}
\bu\bu\bu\bu\sum_{a= (a_1,a_2,\cdots,a_{\l(\mu)})\in DR(\mu)}
\prod_{i\in S(a)}[i]_q.
\cr
}
$$
This completes the proof of our proposition and the proof of  conjecture $\bf (1)$.
\supereject

\noindent{\bol 2. Proof of conjecture (2)}
\sa

Letting ``$\eee$'' denote  the variable which takes the value $-1$ outside the plethystic bracket, 
the modified Hall-Littlewood operator $B_a$
used in the statement of the  compositional 
shuffle conjecture [13] is defined by setting, for any symmetric function $F[X]$,
$$
B_a F[X] = F\big[X+\epsilon\tttt{1-q\over z} \big]\sum_{r\ge 0}z^r e_r[X]\Big|_{z^a}.
\eqno 2.1
$$
In Haglund's book [12] it is shown that  the symmetric  polynomial
$$
B_\mu[X;q]\ses B_{\mu}\, 1 \ses B_{\mu_1}B_{\mu_2}\cdots B_{\mu_{\l(\mu)}}\, 1
$$ 
indexed by any partition  $\mu$ is, up to a factor, the LLT polynomial of a Balanced path indexed by $\mu$.  A Dyck path $D$ is said to be balanced if and only if every North segment of $D$ is immediately followed by an East segment of  equal length.

Our goal in this  section is  to prove the following
\sas

\noindent{\bol Theorem 2.1} 

{\ita For any composition $p$ we have
$$
B_p\, 1\Big|_{q=1+q}\ses \sum_{\mu\part |p|} e_\mu P_\mu(q)
\eqno 2.2
$$
for some polynomial $P_\mu(q)\in \BN[q]$}.
\sa
We must mention that this $e$-positivity is quite surprising since for some compositions, the polynomial $B_p\, 1$ is not even Schur positive.
We will derive this result from the following auxiliary fact.

\sas
\noindent{\bol Proposition 2.1}

{\ita For any integer $a\ge 1$ and partition  
$\mu$ we have
$$
B_a\, e_\mu\Big|_{q=1+q}\ses 
\sum_{\nu\part |\mu|+a}e_\nu\, Q_{\nu,\mu,a}(q)
\eqno 2.3$$
for some polynomial $Q_{\nu,\mu,a}(q)\in \BN[q]$}.
\sas

In fact, since the definition in 2.1 gives 
$
B_a\, 1= e_a
$, we can proceed by induction on the number of components of $p$ and assume 
that 2.2 is valid with $P_\mu(q)\in \BN[q]$.   
Thus we may write
$$
B_p\, 1=\sum_\mu e_\mu P_\mu(q-1).
$$ 
Then for any integer $a$ we derive that
$$
B_a \, B_p\, 1 =\sum_\mu B_a e_\mu P_\mu(q-1),
$$
and 2.3 gives
$$
B_a \, B_p\, 1 =\sum_\mu   \Big(\sum_{\nu\part |\mu|+a}e_\nu\, 
Q_{\nu,\mu,a}(q-1)\Big) P_\mu(q-1)= \sum_{\nu\part |\mu|+a} e_\nu
\sum_\mu Q_{\nu,\mu,a}(q-1)  P_\mu(q-1).
\eqno 2.4$$
Since $Q_{\nu,\mu,a}(q )  P_\mu(q )\in  \BN[q]$
for all $\mu,\nu$  and $a$, the identity in 2.4 completes the induction. This shows  
that we only need to prove Proposition 2.1.
\sas

To this end let us recall that for any expression  $E$ we have
$$
s_\la[X+E]\ses \sum_{\mu \con \la}s_{\la/\mu}[X]s_\mu[E].
$$
In the case that $E=y$ (a monomial) or $E=-y$ 
we obtain
$$
s_\la[X+y]
\ses\sum_{k\ge 0}s_{\la/[k]}[X]\,y^k\ses \sum_{k\ge 0}h_k^\perp s_{\la}[X]\,y^k
$$
and
$$
s_\la[X-y]\ses \sum_{k\ge 0}
s_{\la/[1^k]}[X]\, (-y)^k\ses \sum_{k\ge 0}e_k^\perp s_{\la}[X](-y)^k.
$$
and since the Schur functions are a basis, for any symmetric function $F[X]$  we can write
$$
F[X+y]=\sum_{k\ge 0}y^k h_k^\perp F[X]
\ess\ess\ess\ess\ess\ess and \ess\ess\ess\ess\ess\ess
F[X-y]=\sum_{k\ge 0}(-y)^k e_k^\perp F[X].
$$
Therefore, we derive 
$$
F\big[X+\eee\tttt{1-q\over z} \big]=
F\big[X+\tttt{\eee\over z}-\tttt{\eee q\over z} \big]\ses \sum_{r,s\ge 0}(-1/z)^s(q/z)^r
e_r^\perp h_s^\perp F[X].
$$
The operator in 2.1 can then be rewritten as
$$
B_a\ses \sum_{r,s\ge 0}(-1)^s q^r e_{a+r+s} e_r^\perp h_s^\perp.
\eqno 2.5
$$
To compute $B_a e_\mu$ we will depict $e_\mu$ as the skew Schur function obtained by juxtaposing, corner to corner and on top of each other, columns of lengths $\mu_1,\mu_2,\ldots,\mu_{\l(\mu)}$. For instance the $e$-basis element $e_3e_2e_1$ will be depicted as the leftmost skew diagram in  the following display. Given $r\ge 0$ and $s\ge 0$, we now construct a set  $T^{r,s}_\mu$ of labeled tableaux of shape $\mu$ as drawn on the left of the following display. Each element $S\in T^{r,s}_\mu$ has a weight $wt(S)$ which will give
$$
B_a e_\mu\Big|_{q=1+q}  = \sum_{0\le r+s \leq |\mu|} e_{a+r+s} \sum_{S \in T^{r,s}_\mu} wt(S).
$$
To construct $T^{r,s}_\mu$, first select $s$ cells  which are on the top of their columns and inscribe the cells with ``$-1$''. For instance, if $s=2$, we have the following three choices for filling the example on the left:
$$
 \includegraphics[height=1 in]{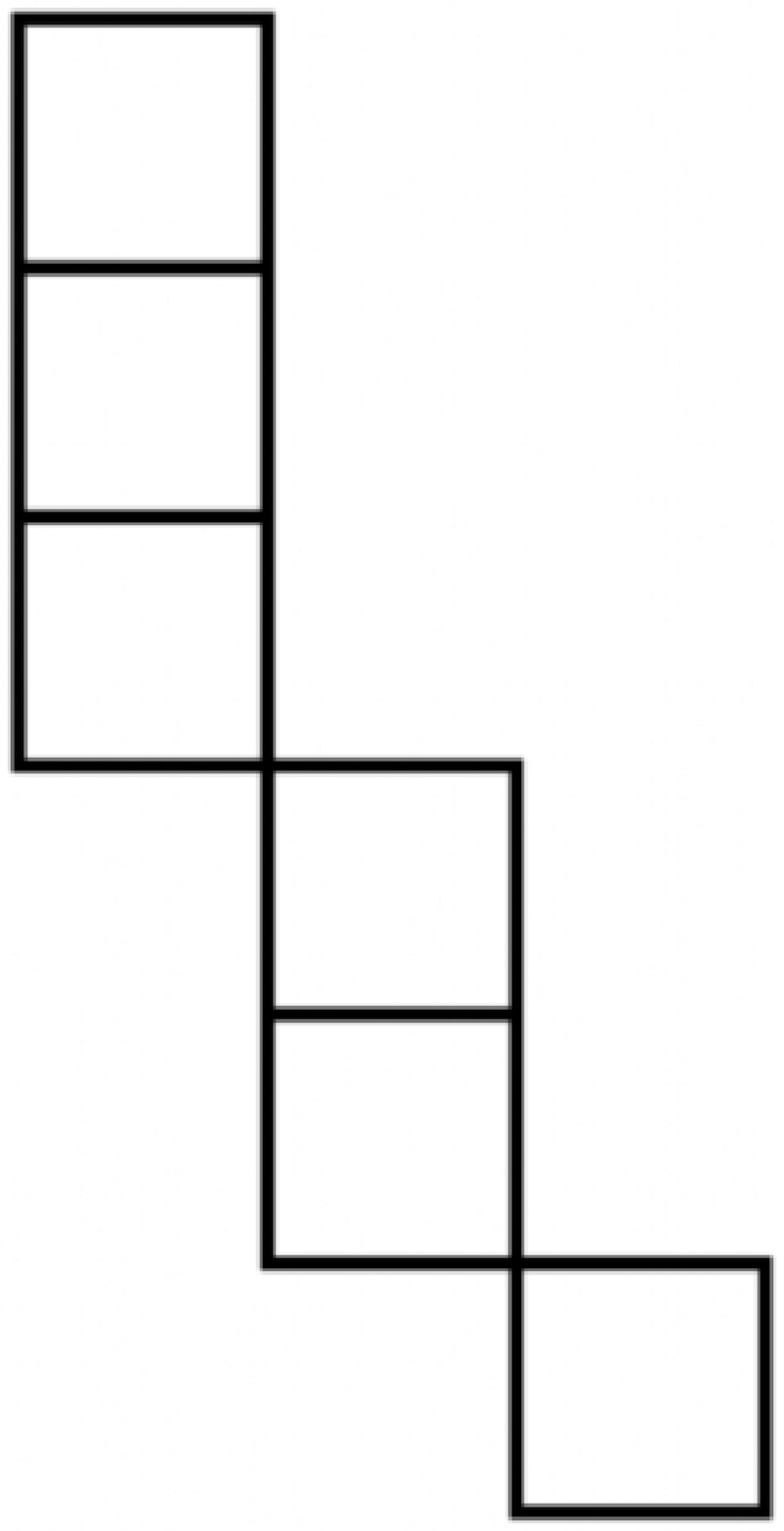}
 \bigsp  \ess\ess\ess   
 \includegraphics[height=1 in]{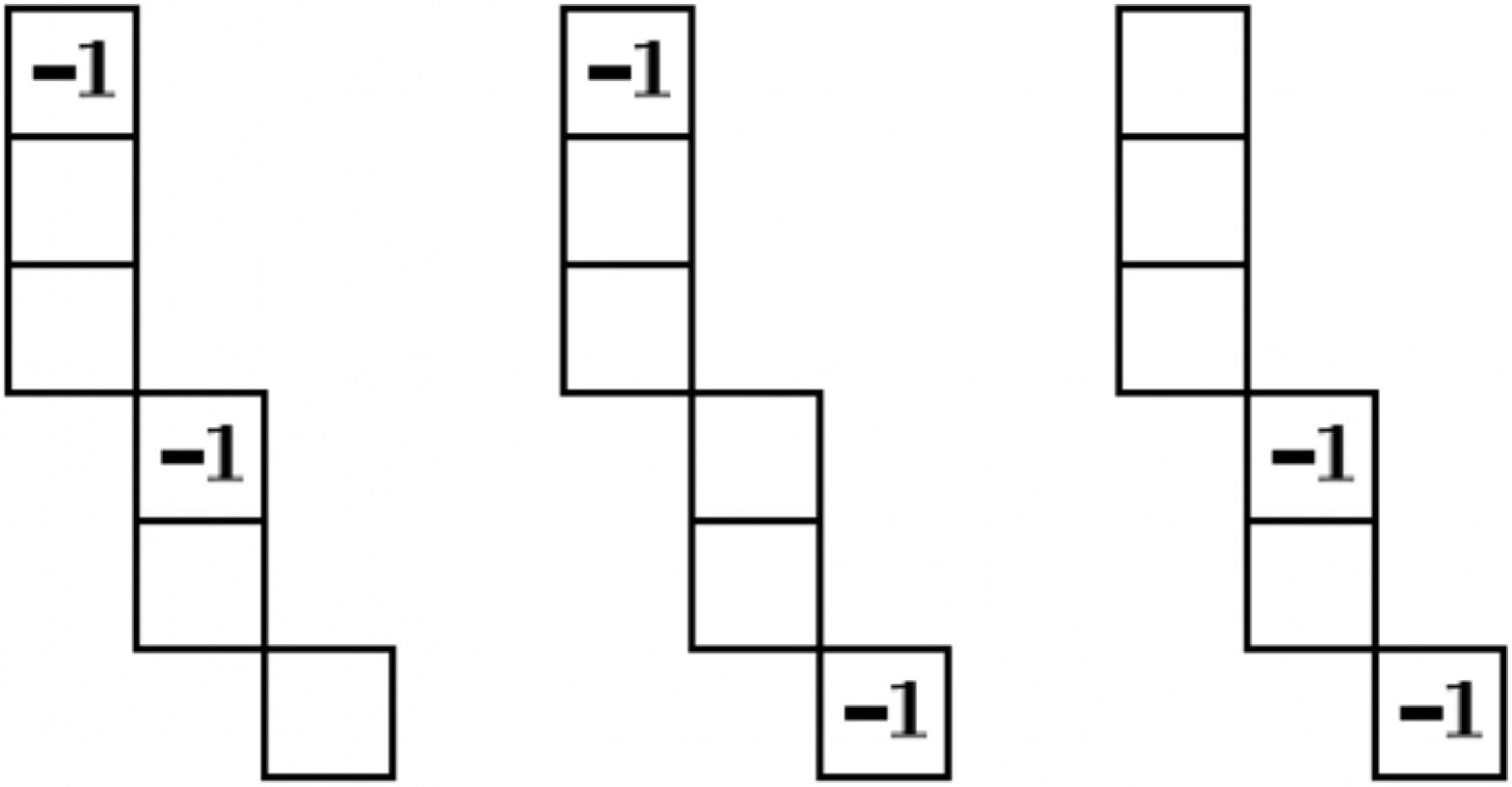}
\ \bigsp  \ess\ess\ess 
 \includegraphics[height=1 in]{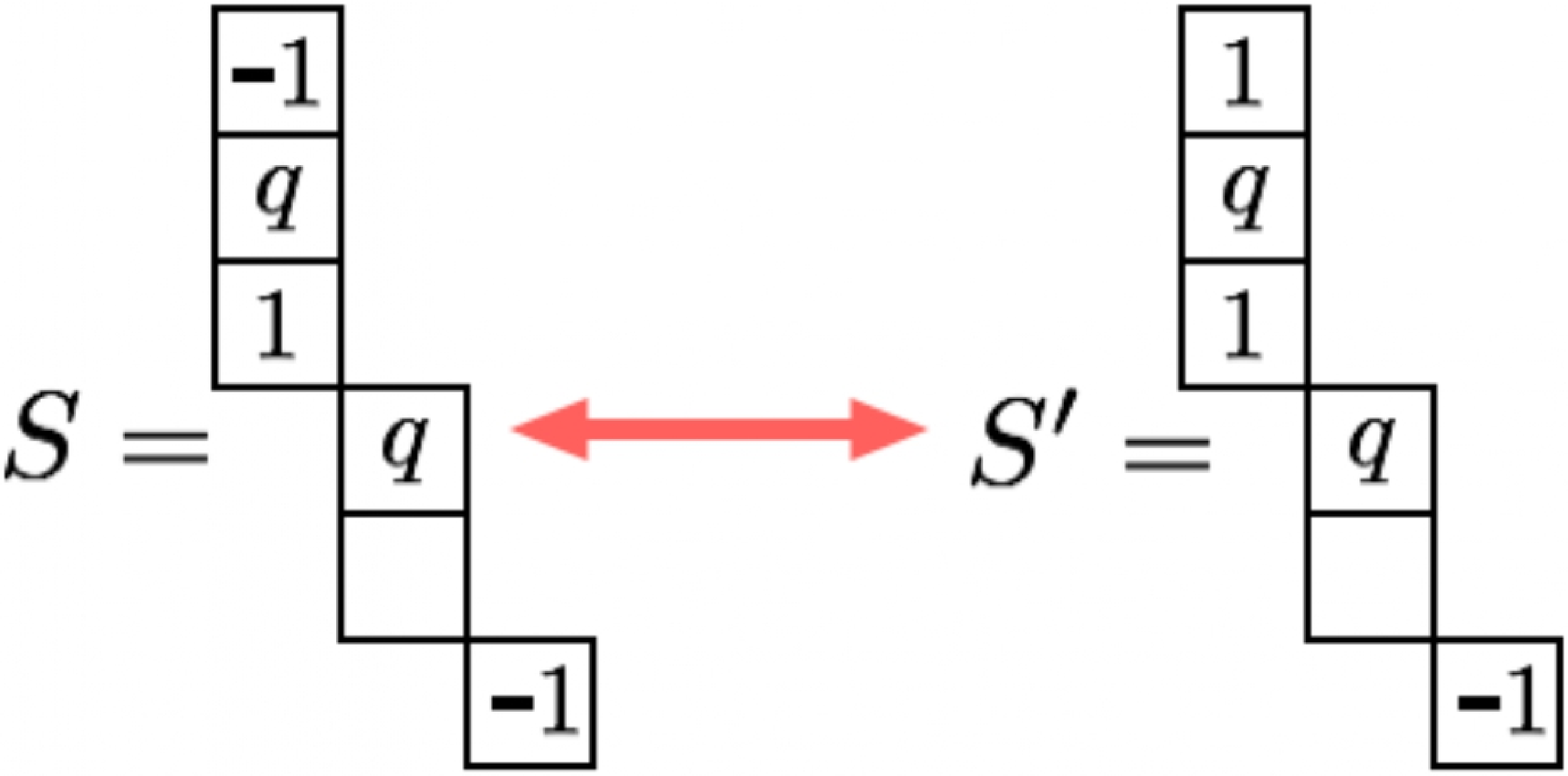}
$$
Next choose $r$ cells so that they form a skew column in the remaining shape, and
for each cell choose  whether to inscribe it with a ``$1$'' or ``$q$''. One example with $s =2$ and $r=3$ is given by
the left member of the rightmost pair. Let $\lambda(S)$ be the partition whose parts are the numbers of empty cells in the columns of $S$. The above example would then produce the partition $(1)$ since there is one empty cell in column $2$. Let $|S|$ be the product of the entries in the cells of $S$. 
The weight of this object is computed by taking
$$
wt(S) = |S| \cdot e_{\lambda(S)}.
$$
The example above would give $wt(S) = (-1) \cdot q \cdot 1 \cdot q \cdot (-1) \cdot e_1 = q^2 e_1$ and 
$wt(S') = (-1) \cdot q \cdot 1 \cdot q \cdot 1 \cdot e_1 = -q^2 e_1.$
\supereject

We now show that $\sum_{S \in T^{r,s}_\mu} wt(S)$ is a positive polynomial by a sign-reversing involution. 
Given $S$, scan from left to right for the first top cell in a column that is either inscribed with a $1$ or a $-1$. Switch the $1$ into a $-1$ in the first case, and switch the $-1$ to a $1$ in the second case. If no such entry exists, leave the tableaux fixed. This is clearly an involution, and it is sign-reversing since we are negating the value of $|S|$, yet preserving the number of $q$'s. 
This involution pairs off the two labeled diagrams 
in the above display.

 Let $U^{r,s}_\mu$ be the subset of $T^{r,s}_\mu$ with the condition that if the top cell of a column is labeled, then it contains a $q$. Thus we have 
 $$
B_a e_\mu \Big|_{q=1+q} = \sum_{r+s \leq |\mu|} e_{a+r+s} \sum_{S \in U^{r,s}_\mu} wt(S),
$$
which is a positive polynomial, completing our proof of Proposition 2.1.
\sa

Theorem 2.1 has a beautiful application. To state it we need some auxiliary facts and notation.
\sas

In a recent posting in the ArXiv Mike  Zabrocki [18] states a general conjecture   asserting that a certain $S_n$ module has the symmetric function appearing in the Delta conjecture [14] as Frobenius Characteristic.
Our application  of  Theorem 2.1 is  that
Zabrocki's Conjecture implies that 
a submodule of Zabrocki's module exhibits 
the $e$-positivity phenomenon.
\sas

  To define Zabrocki's submodule  we consider the vector space  
$R_n[X,\Theta]=\BQ[x_1,\ldots,x_n;\theta_1,\ldots,\theta_n]$, with the $x_i$ commuting variables, the $\theta_j$ anti-commuting and commuting with the  $x_i$. This space is itself an $S_n$  module under the diagonal action. The latter is simply defined by  letting a permutation $\sig\in S_n$ send $x_i$ into $x_{\sig_i}$ and $\theta_i$  into 
 $\theta_{\sig_i}$. The Zabrocki submodule is none other than the analogue of the Diagonal Harmonics module when $R_n[X;Y]=\BQ[x_1,\ldots,x_n;y_1,\ldots,y_n]$  
with the $x_i,y_j$ commuting variables.
This given,  Zabrocki's submodule is the quotient of $R_n[X;\Theta]$ by the ideal $I_n$ generated by the  diagonal invariants in $R_n[X,\Theta]^{Sn}$ with vanishing constant term. Let us call them the $x,\theta$-Coinvariants. Or equivalently, by taking the the orthogonal complement of $I_n$, the $x,\theta$-Diagonal Harmonic's.  Now it follows from Zabrocki's conjecture that the Frobenius Characteristic 
of this module is the symmetric function
$$
DH_{x,\theta}[X;q]\ses\sum_{k=1}^n(-t/q)^{n-k}
\om E_{n,k}[X;1/q].
\eqno 2.6
$$
The polynomials $E_{n,k}[X;q]$ were originally defined in [6] by the Pochhammer expansion
$$
e_n\big[X\tttt{1-z\over 1-q}  \big]
\ses \sum_{k=1}^n\tttt{ (z;q)_k \over (q;q)_k}
E_{n,k}[X;q].
\eqno 2.7
$$
In [13]  it is  shown that
$$
E_{n,k}[X,q]\ses \sum_{\aaa\models n} C_{\aa_1}C_{\aa_2}\cdots C_{\aa_k}\, 1,
\bigsp(\hbox{$\l(\aaa)=k$}),
\eqno 2.8
$$
where the operators $C_a$ are defined by setting
 for any symmetric function 
$F[X]$,
$$ 
C_a F[X]\ses (-{\tttt{ 1\over q }})^{a-1} 
F\big[X-\tttt{1-1/q \over z}
\big]\sum_{r\ge 0} z^r h_r[X]\Big|_{z^a}.
\eqno 2.9
$$
It turns out  that we can derive an expression   similar to 2.8 for the polynomial
$$
\widetilde E_{n,k}[X,q]\ses (-1/q)^{n-k}\om E_{n,k}[X,1/q].
\eqno 2.10
$$
To see this, notice first that the operators
 $B_a$ defined in 2.1 may also be defined by setting $B_a=\om \widetilde  B_a\om$ with
$$
\widetilde B_a  F[X]\ses F\big[ X-  \tttt{1-q \over z}\big]]\sum_{r\ge 0} z^r h_r[X]\, \Big|_{z^a}.
\eqno 2.11
$$
This given, notice that replacing $q$ by $1/q$ in 2.9 we obtain
$$ 
C_a^{1/q} F[X]\, =\, (- q )^{a-1} 
F\big[X-\tttt{1- q \over z}
\big]\sum_{r\ge 0} z^r h_r[X]\Big|_{z^a} =\,
(- q )^{a-1}\WB_a F[X]
,
\eqno 2.12
$$
\vskip -.2 in
\noindent
and 2.8 becomes
$$
E_{n,k}[X,1/q]\ses (- q )^{n-k}\sum_{\aaa\models n} \WB_{\aa_1}\WB_{\aa_2}\cdots \WB_{\aa_k}\, 1.
\eqno 2.13
$$
Since $\om 1=1$, the equality in 2.13 can also be rewritten as
$$
\om E_{n,k}[X,1/q]\ses (- q )^{n-k}\sum_{\aaa\models n} \om\WB_{\aa_1}\om\, \om\WB_{\aa_2}\om\cdots \om\WB_{\aa_k}\om \, 1,
\eqno 2.14
$$
\vskip -.2 in
\noindent
and 2.10 becomes 
$$
\widetilde E_{n,k}[X,q]\ses
 \sum_{\aaa\models n} B_{\aa_1}B_{\aa_2}\cdots B_{\aa_k}\, 1.
\eqno 2.15
$$
Thus the Zabrocki  conjecture in 2.6
may be also computed by the formula
$$
DH_{x,\theta}[X;q]\ses \sum_{k=1}^n t^{n-k}
\sum_{\aaa\models n} B_{\aa_1}B_{\aa_2}\cdots B_{\aa_k}\, 1.
\eqno 2.16
$$
Since Theorem 2.1 states that all the symmetric polynomials
$\ess
 B_{\aa_1}B_{\aa_2}\cdots
B_{\aa_k}\, 1
\ess $
exhibit  the $e$-positivity phenomenon, it follows from 2.16 that 
so does the polynomial 
$
x\theta DH[X;q,t].
$

Since computer data shows that the summands  in 2.16 are not necessarily Schur positive when $\aa$ is not a partition, the question remains as to what mechanism causes their sum 
to be Schur positive. The answer is quite simple. In fact, let us recall that in Haglund's book [12] it is shown that the polynomial
$B_{\aa_1}B_{\aa_2}\cdots B_{\aa_k}\, 1$ is the LLT polynomial of a Dyck path only when the $\aa_i's$ are the components of a partition. However, it follows by combining Zabrocki's conjecture in [18], the conjecture in 2.6 and the Delta conjecture at  $t=0$ (see [7]) that the left-hand side   
2.6 can be  rewritten in the two forms 
$$
DH_{x,\theta}[X;q]\, =\,\sum_{k=1}^n z^{n-k}\sum_{a \models n}B_{\aa_1}B_{\aa_2}\cdots B_{\aa_k}\, 1\,=\, \sum_{D\in \CD_n}t^{area(D)}LLT_D[X;q]Hag(D;z)\Big|_{t=0} 
\eqno 2.17
$$
where $Hag(D;z)$ is the Haglund factor of $D$. If the component   $D_i$ is  the co-area  of the $i^{th}$ North step of $D$  and $a_i$ denotes the number of area cells in the $i^{th}$ row, then we can write
$$
Hag(D;z)= \prod_{i=2}^n (1+z/t^{a_i})^{\chi(D_i=D_{i-1})},\ess\ess  
\hbox{ thus}\ess\ess\ess 
t^{area(D)}Hag(D;z)=\prod_{i=2}^n  t^{a_i\chi(D_i\neq D_{i-1})} \prod_{i=2}^n (t^{a_i}+z)^{\chi(D_i=D_{i-1})}.
$$  
In order that the latter factor contributes to the sum in 2.17 it must be that whenever $D_i\neq D_{i-1}$ then $a_i=0$. That means the path $D$ must hit the diagonal every time $D_i\neq D_{i-1}$. This forces $D$ to be a Balanced Dyck path and the equality in 2.17 to be none other than
$$
DH_{x,\theta}[X;q]\, =\,\sum_{k=1}^n z^{n-k}\sum_{a \models n 
}B_{\aa_1}B_{\aa_2}\cdots B_{\aa_k}\, 1\,=\,\bu \sum_{D\in \CD_n\,\, balanced}\bu\bu LLT_D[X;q]\,\,z^{\sum_{i=2 }^n\chi(D_i=D_{i-1})}.
\eqno 2.18
$$
Since $\sum_{i=2 }^n\chi(D_i=D_{i-1})=n-k$ if $D$ hits the diagonal in $k$ places, it follows that we also have the equality
$$
\sum_{a \models n}B_{\aa_1}B_{\aa_2}\cdots B_{\aa_k}\, 1\,=\,\bu \sum_{D\in \CD_n\,\, balanced \,\, with\,\, k\,\, hits }\bu\bu\bu\bu LLT_D[X;q].
\eqno 2.19
$$ 
This clearly explains why the left hand side ends up being Schur positive.
\supereject
\noindent{\bol 3. Some consequences and conjectures} 
\sas
 
 First and foremost, the $e$-positivity phenomenon 
  suggests that an action of $S_n$ is   involved, yet so far none of our proofs uses anything of the sort. Let us recall that any finite group action breaks up into a  direct sum of transitive group actions. Moreover, transitive submodules of group actions are none other than the orbits of the action. However orbit actions are equivalent to left coset actions. Thus the character of an orbit action, as an element of the acting group algebra,  can be  simply expressed  in terms of the stabilizer of 
any element  of the orbit. In the case of $S_n$, 
if all these stabilizers happen to be Young subgroups, then the Frobenius characteristic of the 
character of the action must be $h$-positive or $e$-positive. 

These ideas suggest a simple computer exploration. Namely,  finding  out what is beinq acted upon. 
We can do this by computing the Hilbert series of the conjectured module and set $t,q=1$. In  those cases where we   obtain a sequence of integers, the resulting data is a good candidate for the encyclopedia of integer sequences. The simplest case is $\nabla e_n$ which is the symmetric function side of the shuffle conjecture, now a theorem [2]. Since the combinatorial side is obtained as a sum of all $n\times n$ Dyck paths LLT's, $\nabla e_n$ itself should also exhibit the $e$-positivity phenomenon. Computer data strongly confirms that it does. Thus it seems worthwhile to find out what is  the cardinality of the set of objects that $S_n$ should be acting upon in this case. Doing this exploration with $\nabla e_n\Big|_{q=2}$ for
 $n=1,2,3,4,5,6,7$ we get   
$$
1,\, 4,\,38,\,728,\,26704,\,1866256,\, 251548592,\, \ldots
$$ 
Entering this sequence in the encyclopedia returns an avalanche of hits. The immediate answer is
\sas

\centerline {\ita The number of connected graphs on $n+1$ nodes.}
\sas

There is even a connection with Novak's Free probability notes [17] where we can  find a list of all the 38 connected graphs on $4$ nodes. A further search more closely connected with the replacement $q\RA 2$ yields the papers of Kreweras [15] and Gessel-Wang [10] who now appear to have hit the tip of an iceberg. 

Finding a bi-graded $S_n$-Module with Frobenius characteristic the $e$-basis expansion of the polynomial $\nabla e_n\big|_{q=1+q}$ would make an interesting research problem indeed.
Likewise, the conjecture that the LLT polynomials of Dyck paths exhibit the $e$-positivity phenomenon suggests that such a module should exist  

\vskip .05in
\hfill 
$ \includegraphics[height=1.1in]{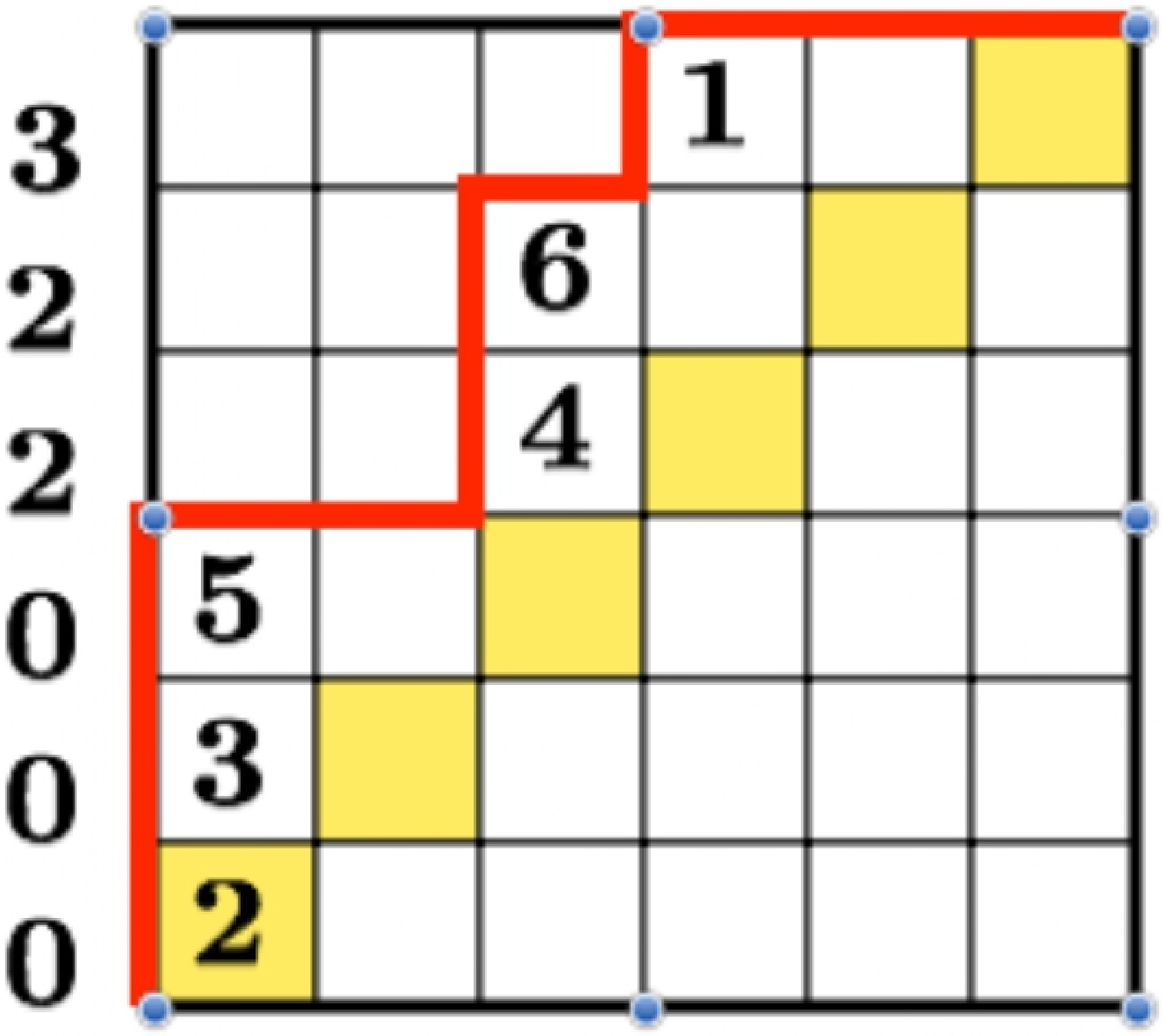}$

\hsize 5.1in
\vskip -1.18in
\noindent
also in these cases. But before we focus more closely on Dyck path LLT's it will be good to recall the  definition of the ingredients that enter in their construction.   
In the adjacent display we have our depiction of a Parking Function. To begin we have drawn a Dyck path $D$ in the $6 \times 6$ lattice square $R_6$. This is a path that goes from $(0,0)$ to $(6,6)$ by unit North and East steps always remaining weakly above the lattice diagonal (the yellow cells). We have also labeled the cells adjacent to the

\vskip -.02in
\hsize 6.5in
\noindent
North steps of $D$ by the numbers $\,\, 1,\, 2,\,3,\,4,\,5,\,6$, usually referred to as ``cars'' in a  column increasing  manner. We have two statistics of a parking function called $area(PF)$ and $dinv(PF)$. The statistic $area(PF)$ is actually the area of $D$ which is the number of  cells between the path and the lattice diagonal. The statistic $dinv(PF)$ is obtained as the total number of ``primary'' and  ``secondary dinvs.   Two  cars in the same diagonal yield a primary dinv if the one on the left is smaller  than the one on the right. 
A secondary dinv is yielded by two cars when the one on the left is on a higher diagonal but
 adjacent to the diagonal of the car on the right, 
and the car on the left is larger than the  car on the right. In the above example we have two primary dinvs $3,4$ and $5,6$ and the secondary dinv   $5,4$. The word of the parking function, denoted $\sig(PF)$, is the permutation obtained by reading the cars by diagonals from right to left starting from the highest and ending with  the lowest. Thus for our example 
$\sig(PF)=165432$. The largest dinv is obtained when $\sig(PF)=654321$. This is  the dinv of     
the Dyck path.

\supereject

This given, the following identity gives us precise information as to the number of orbits (or $e$-basis elements) and their weight. 
\sas
\noindent{\bol Proposition 3.1}

{\ita Suppose for a Dyck path $D$ in the $n\times n $ lattice square we have the expansion
\vskip -.1in
$$
LLT_D[X;1+q]\, = \sum_{\mu\part n}e_\mu[X]P_\mu(q),
\eqno 3.1
$$

\vskip -.18in
\noindent
then}
$$
\sum_{\mu\part n}P_\mu(q)\ses (1+q)^{dinv(D)}.
\eqno 3.2
$$

\vskip -.2in
\noindent{\bol Proof}

Recall that by definition   the LLT polynomial of a Dyck path $D$ is given by the formula

\vskip -.1in
$$
LLT_D[X;q]\ses\sum_{D(PF)=D}q^{dinv(PF)} F_{pides(PF) }[X]
\eqno 3.3
$$

\vskip -.092in
\noindent
where the sum is over all  Parking Functions supported by $D$. The last factor here is the Gessel quasi-symmetric function basis element indexed by $pides(PF)$, the composition giving the descent set of the inverse of the word $\sig(PF)$. Since $LLT$ polynomials are symmetric it follows from a theorem of Gessel [9] and 3.3 that 

\vskip -.155in
$$
\LL LLT_D[X;q]\scs s_{1^n}[X]\RR\ses 
\sum_{D(PF)=D}q^{dinv(PF)}\chi(pides(PF)=1^n),
\eqno 3.4
$$  

\vskip -.09in
\noindent
where the left hand side of this identity is a Hall scalar product of two symmetric functions.
On the right hand side the equality $pides(PF)=1^n$ can only happen when $\sig(PF)=n\cdots 321$. 
This reduces the sum in 3.4  to a single term and in that case we have $dinv(PF)=dinv(D)$.
Thus 3.4 may also be rewritten as 
$$
\LL LLT_D[X;1+q]\scs s_{1^n}[X]\RR\ses (1+q)^{dinv(D)}.
$$
But then 3.2 follows from 3.1 since for any $\mu\part n$ we have 
$
\ess\LL\, e_\mu[X] \scs s_{1^n}[X]\, \RR= 1.
$ 
This completes our proof.
\sa

Guided by the identity in 3.2 and supported by   computer data we have been led to the following
\sas

\noindent{\bol Conjecture 3.1}

{ Given a Dyck path $D$ in the $n\times n$ lattice square $R_n$, the following algorithm 
constructs the $e$-basis expansion of the polynomial $LLT_D[X;1+q]$.  Draw the parking function $PF$ with $\sig(PF)=n\cdots 321$ and determine the set of pairs of cars $(a,b)$ with $a<b$ producing a dinv. Call this ``$dinvset(D)$''.
Constructs the set of pairs of cars $a<b$ placed one above the other in one of the columns of $D$ and call it 
``$forced(D)$''. Then our final product can be written  in the form

\vskip -.2in
$$
LLT_D[X;1+q]=\sum_{S\con dinvset(D)}q^{|S|}
e_{\mu(S)}.
\eqno 3.5
$$

\vskip -.1in
\noindent
To construct the terms of this sum repeat 
 the following $4$ steps until all the subsets of $dinvset(D)$ have been processed. Begin by setting $out:=0$.
\sas

\item{\bol (1)} {\ita  Choose a subset $S$ of $dinvset(D)$ and set $temp:=q^{|S|}$.
\sas
\item {\bol (2)} {\ita Using all the pairs in $S$ and $forced(D)$ to construct the poset $\Pi=\big(\{1,2,\ldots, n\},\preceq\big)$ (here each pair $(a,b)\in S\cup forced(D)$ must be interpreted as $a\prec b$)}.

\item {\bol (3)} Recursively start by setting
$\Pi':=\Pi$ and $max(\Pi'):=n$ and repeat the following commands until $\Pi'$ has no more }elements:
\itemitem { (i)} {\ita Determine the downset of $max(\Pi')$, (the set of elements of $\Pi'$ that are $\preceq max(\Pi')$)
\itemitem { (ii)} if the size of this downset is $k$ do $temp:=temp\times e_k$, 
\itemitem { (iii)} Remove from $\Pi'$ all the elements of $downset\big(max(\Pi')\big)$ and let the result be the new $\Pi'$,
\sas
\item {\bol (4)} Save the result of   steps 
$\bf (1),(2),(3)$ by the command $out:=out+temp$.}
\sas

\noindent
We conjecture that at the completion of this algorithm $out$ will give the right hand side of 3.5. 
\supereject

\sas

\hfill 
$ \includegraphics[height=1.2 in]{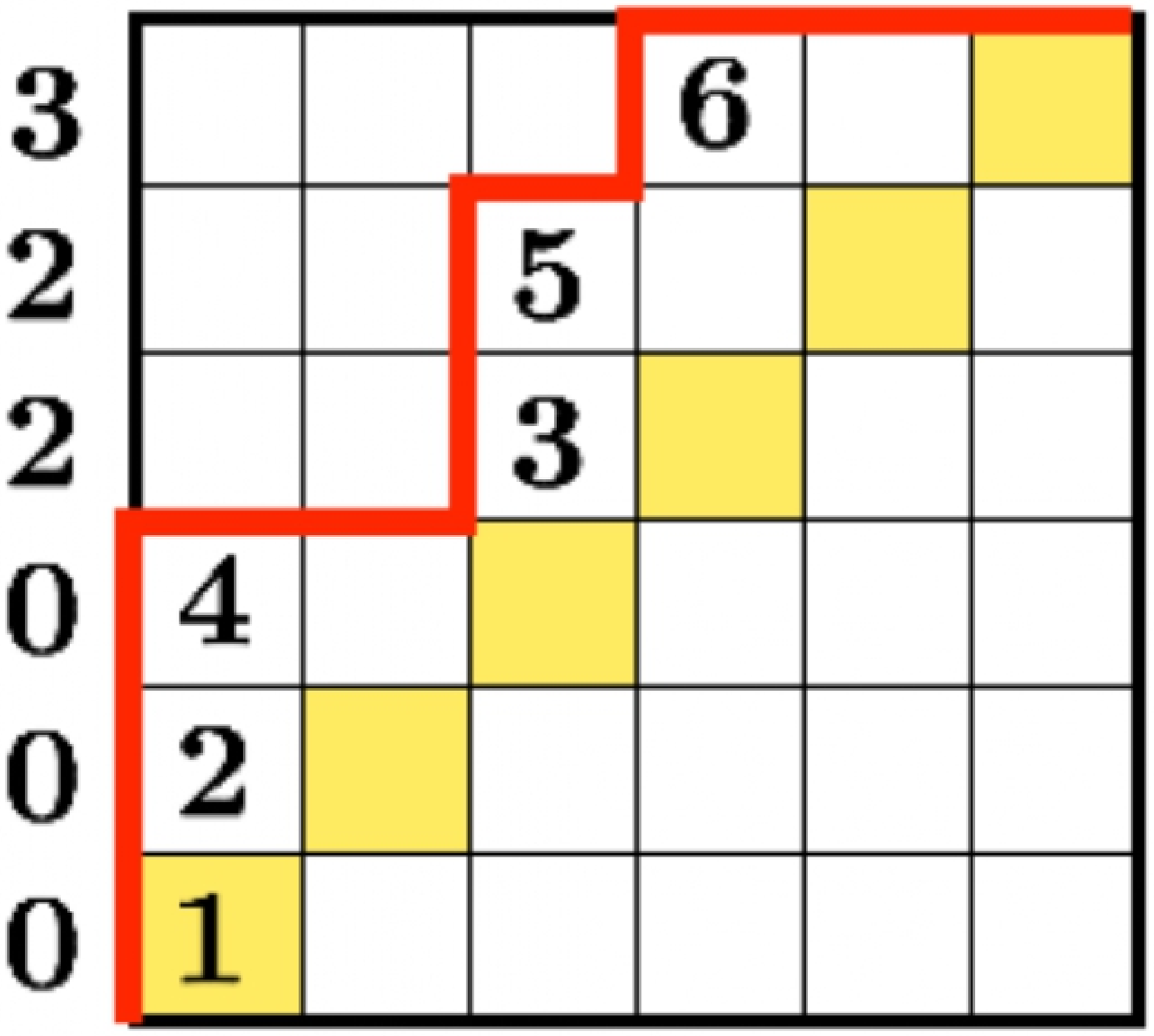}$

\hsize 5in
\vskip -1.3in

This is best illustrated by an example.
 In the  display on the right we have depicted in $R_6$  the same Dyck path $D$ we had in our previous display.  However we have labeled the cells adjacent to the North steps of $D$ by the permutation  $654321$. More precisely, in diagonal $2$ we placed $654$, in diagonal $1$ we placed 
$32$   and in diagonal $0$ we placed $1$.  This labeling gives the Parking Function with the highest dinv.   It causes a total of $5$ dinvs which include primary  pairs $(2,3)$, $(4,5)$, $(4,6)$, $(5,6)$ and a secondary  pair $(3,4)$. 
Thus  we obtain

\vskip -.2in
\hsize 6.5in
$$
a)\ess\ess dinvset(D)=\big\{(2,3),(4,5),(4,6),(5,6),(3,4)\big\}
\ess\ess \hbox{and}\ess\ess
b)\ess\ess forced(D)=\big\{(1,2),(2,4),(3,5)\big\}.
\eqno 3.6
$$
Since the number of subsets of $dinvset(D)$ is $ 2^5$ we will only carry out  steps $\bf (1),(2),(3)$ of the above algorithm in the special cases
$S=\{(4,6)\}$ ,  $S=\{(2,3),(4,5),(4,6)\}$, $S=\{(3,4), (4,6)\}$ ,$S=\{(3,4)\}$ and obtain the following four posets.

\vskip -.30in
$$
 \includegraphics[height=1.1 in]{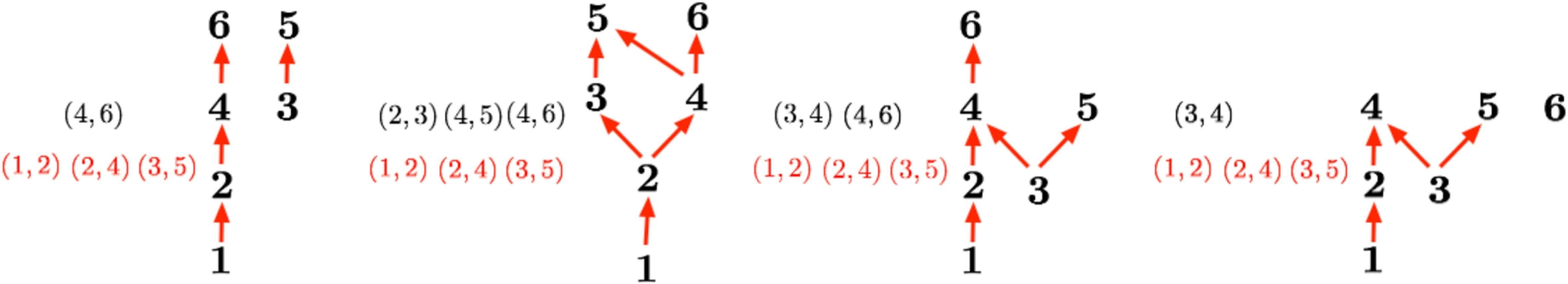}
$$

\vskip -.2in 
 \noindent
Processing these posets gives the following four terms of the polynomial in 3.5 for our choice of $D$:
$$
q\,  e_4e_2,\ess\ess\ess\ess
q^3 e_4e_2,\ess\ess\ess\ess
q^2e_5e_1,\ess\ess\ess\ess
q\, e_3e_2e_1.
 \eqno 3.7
 $$
Using the definition in 3.3 and replacing each Gessel Fundamental by the Schur function indexed by the same composition yields the $e$-basis expansion  (see [3] and [8] for this operation)
$$
\eqalign{
LLT_D[X;1+q]=
(q^5+4q^4+&5q^3+2q^2)e_6+(q^4+4q^3+4q^2+q)e_5e_1+
\cr
&
(q^3+2q^2+q)e_4e_2+
(q^2+q)e_4e_1^2+(q^2+q)e_3^2+
(1+q)e_3e_2e_1
\cr
}
$$
and we can verify the terms in 3.7 do occur in the $e$-basis expansion of our polynomial.  
\sas

\vskip -.2in
\hfill
$
 \includegraphics[height=1.1 in]{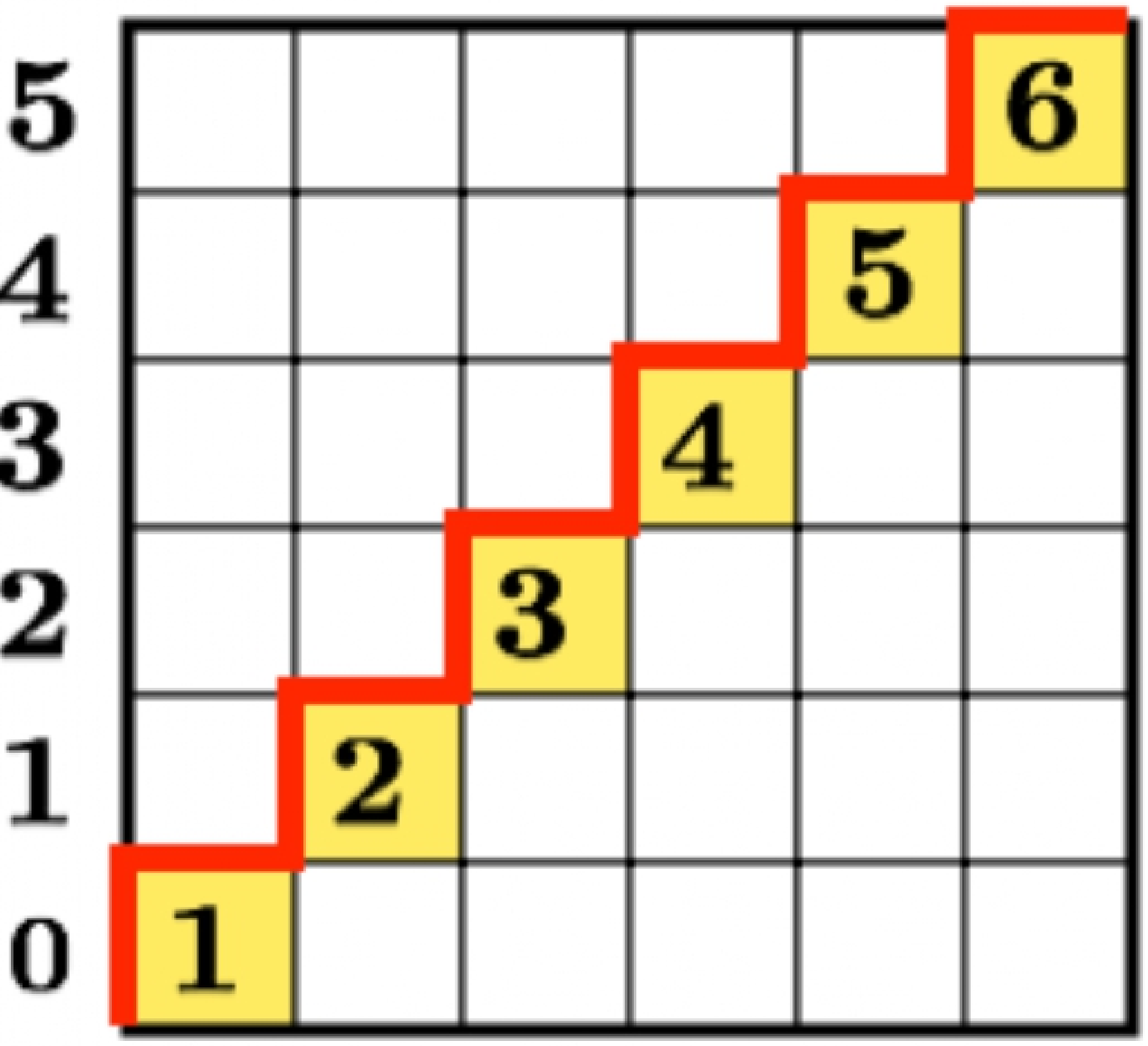}
$

\vskip -.98in
\hsize 5.2 in
It turns out that this  conjecture has a surprising consequence that can actually be proved. This occurs for the Dyck path
$D$ with no area as illustrated in the adjacent display. It easily seen that in this case we have no forced pairs and

\vskip -.1 in
$$
dinvset(D)=\{ \, (a,b)\,:\, 1\le a<b\le n\,\}. 
\eqno 3.8
$$

\hsize 6.5in
\vskip -.03in
\noindent
Let us see what are the possible selections of subsets $S$ of $dinvset(D)$ yielding    
\vskip -.15in
$$
downset(6)=\{t_1<t_2<t_3<6\}.
\eqno 3.9
$$ 
Recall  that in this case the first factor of the final $e$-basis element  is $e_4$. In order that we will be  able to recognize  other factors it will be necessary to work with  the original polynomial $LLT_D[X;q]$. Thus rewriting  $q=q-1+1$ we can interpret $``q-1"$ as including  a given pair in $S$ and $``1"$ as not including it. So the question is how can we guarantee that our choices will result in $downset(6)$  to be  as  in 3.9. We also must recall  that each 
new  edge $(a,b)$ we add to $S$ must satisfy 
$(a<b)$. Thus to assure that $t_{3}$ is in 
$downset(6)$ we need only add to  $S$ the edge $(t_{3},6)$. Recursively, to assure that
$t_2$ is in  $downset(6)$ we need only add  to $S$ at least one of the edges $(t_2,t_3)$ or $(t_2,6)$.  
Finally to assure that  $t_1$ is in $downset(6)$ we need at least one of the three edges $(t_1,t_2)$, $(t_1,t_3)$, and $(t_1,6)$.  

This given, the contribution 
to the $q$ factor of the final polynomial resulting from these three choices is   
$(q-1)(q^2-1)(q^3-1)$. The $-1$ in each case results from the fact that in each case we are required to pick at least one of the given choices,
i. e. picking none is not permitted.
\supereject

But we are not done yet  with powers of $q$. Set
$$
R =\{1,2,3,4,5,6\}-T 
$$
with $T=downset(6) $. Without affecting  $downset(6)$ 
we can add to $S$ any edges $(t,r)$ with $t\in T$ and $r\in R$ with $t<r$. For instance if 
$\{t_1,t_2,t_3\}=\{2,4,5\}$ then $R=\{1,3\}$. The insertion or not insertion of such an edge results in an additional  $q$ for each of the available choices. For this particular instance of $t_1,t_2,t_3$ we have only one choice $(2,3)$. Picking or not picking this pair yields a $q$ factor. In the display below we have listed all possible choices of $\{t_1,t_2,t_3\}$ (vertically), their remainder $R$   and the power of $q$ they contribute.
$$
 \includegraphics[height=.9 in]{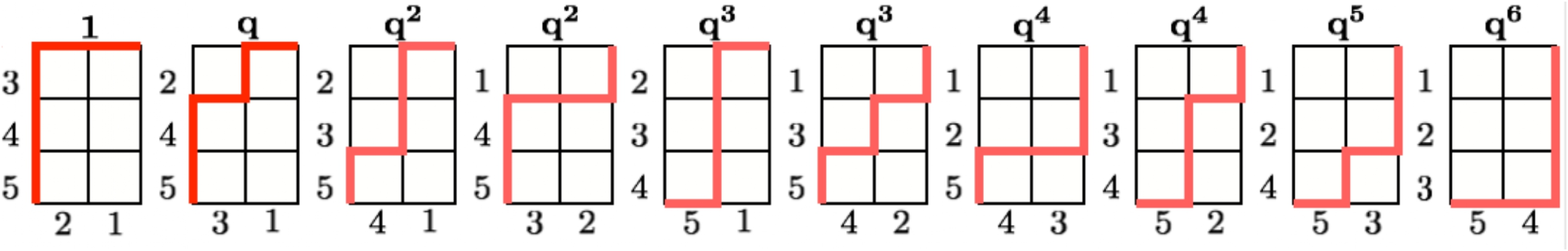}
$$
We can thus easily see that the $q$-factor  that
accounts for all  the  choices that  do  not affect the size of $downset(6)$, is none other than the polynomial that $q$-counts by area the partitions that are contained in a $3\times 2$  rectangle, that is the $q$-binomial coefficient
$$
\Big[{6-1\atop 4-1}\Big]_q.
\eqno 3.10
$$
To complete the contribution due to this $downset$ we need to observe that when $D$ is the  no area path in the lattice square $R_n$ we may as well use the notation
$$ 
LLT_D[X;q]\ses LLT_n[X;q].
\eqno 3.11
$$
The idea is that the choices we were forced to  make to  assure our particular
$downset$ will not affect the remaining construction. More precisely, at this point it is natural to assume that whatever must be added to complete the contribution of this $downset$ will be recursively provided by the construction of the polynomial $LLT_{6-4}[X;q]$.

In summary, the contribution to the polynomial 
$LLT_6[X;q]$ due to all $downsets$ of size $4$ should be
$$
e_4[X]\,(q-1)(q^2-1)(q^3-1)
\Big[{6-1\atop 4-1}\Big]_qLLT_{6-4}[X;q].
\eqno 3.12
$$

It turns out that the validity of the idea that suggested 3.12 is confirmed by
the following recursion we can actually prove.
\sas

\noindent{\bol Theorem 3.1}

{\ita For any $n\ge 2$ we have}
$$
LLT_n[X;q]\ses \sum_{k=1}^ne_k[X]
(-1)^{k-1}(q;q)_{k-1}\Big[{n-1\atop k-1}\Big]_q\, LLT_{n-k}[X;q].
\eqno 3.13
$$
{\bol Proof}

To begin we must observe that in this particular case the   LLT  polynomials  have explicit expressions.  In fact, a moment's reflection reveals that for each $n\ge 1$ we have
$$
LLT_n[X;q]= \TH_{[n]}[X;q,1]\ses
(q;q)_n\, h_n\big[\tttt{X\over 1-q}\big].
\eqno 3.14
$$
Using this identity in 3.13 the definition of the $q$-binomial coefficient gives
$$
(q;q)_n\, h_n\big[\tttt{X\over 1-q}\big]\ses \sum_{k=1}^ne_k[X]
(-1)^{k-1}(q;q)_{k-1}
{(q;q)_{n-1}\over (q;q)_{k-1}
(q;q)_{n-k}}
\, (q;q)_{n-k}\, h_{n-k}\big[\tttt{X\over 1-q}\big] .
$$
Carrying out all the  obvious cancellations we get
$$
(1-q^n)\, h_n\big[\tttt{X\over 1-q}\big]= -\sum_{k=1}^ne_k[X]
(-1)^{k}
 \, h_{n-k}\big[\tttt{X\over 1-q}\big].
$$
However this is the same as 
$$
(1-q^n)\, h_n\big[\tttt{X\over 1-q}\big]= -\sum_{k=0}^n h_k[-X]
 \, h_{n-k}\big[\tttt{X\over 1-q}\big]  
 \sps h_{n}\big[\tttt{X\over 1-q}\big].
 $$
Or better
$$
 -q^n \, h_n\big[\tttt{X\over 1-q}\big]= -\sum_{k=0}^n\, 
 h_{k}\big[\tttt{-X+qX\over 1-q}\big] 
 h_{n-k}\big[\tttt{X\over 1-q}\big] \,
 = - \, h_n\big[\tttt{qX\over 1-q}\big]. 
 $$
Completing our proof.
\sas

 In trying to prove the recursion of Theorem 3.1 for general LLT's it was discovered that some of the terms recursively constructed were column LLT's. These findings  resulted in further  discoveries.
 \sas

\centerline{\ita  a) The ``areaprime" construction of  column LLT's} 
\sas

\hfill
$
 \includegraphics[height=1.4 in]{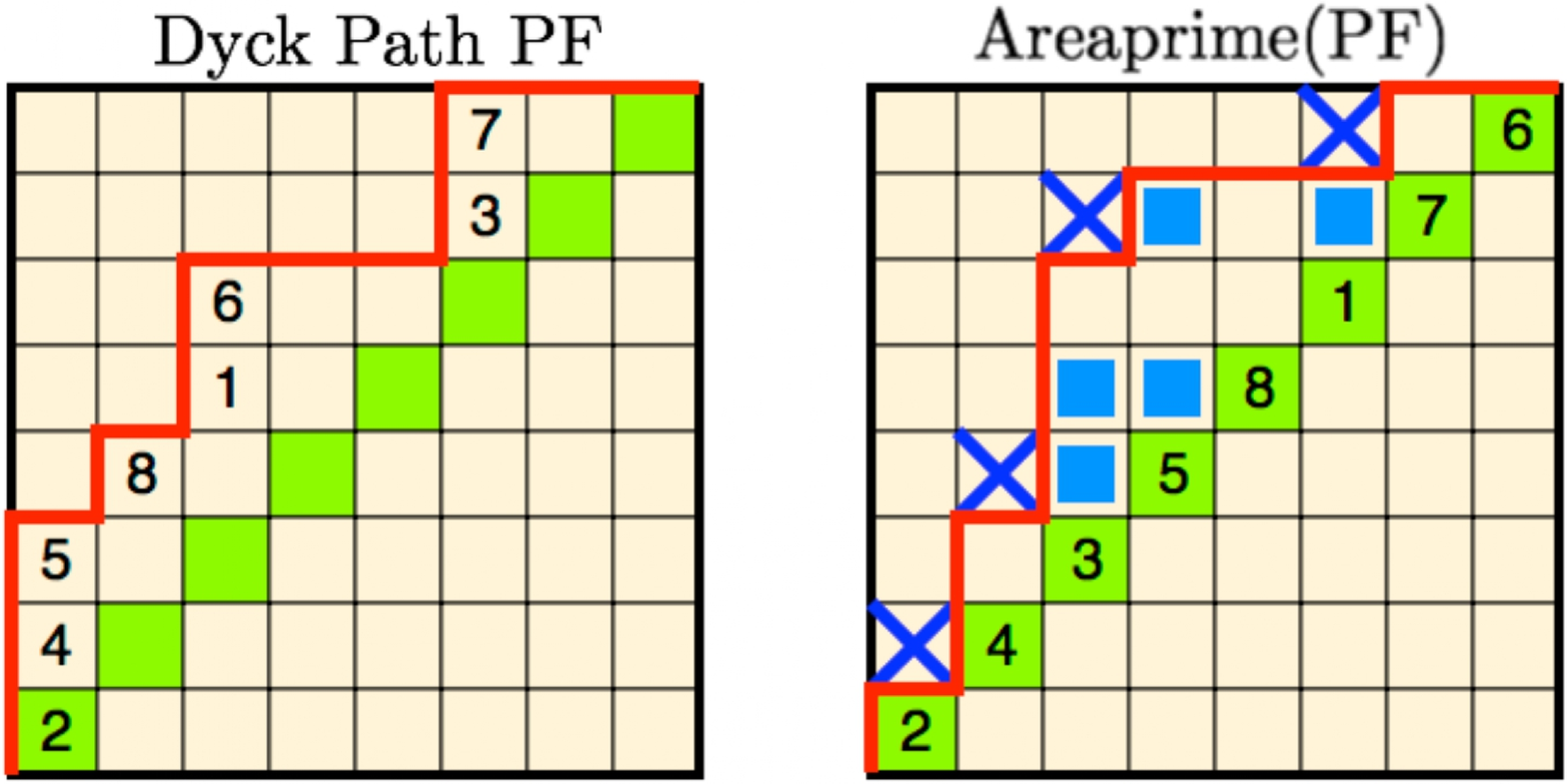}
$

\vskip -1.4 in
\hsize 3.5in

We will start by recalling the construction of the ``areaprime" image of a Parking function. In the adjacent display we have a parking function supported by a Dyck path $D$ and 
its areaprime image. To  obtain this image the first step is to construct the permutation on the diagonal. This is done by reading the cars in $PF$ by diagonals from left to right   
starting from the lowest and ending with 
the

\hsize 6.5in
\noindent
 highest. This done, we 
determine the positions of the Blue crosses. Each cross is determined by a pair $(a,b)$  of cars $a<b$ with $a$ directly below $b$ in one of the North segments of $D$. These pairs are  $(2,4),(4,5),(1,6), (3,7)$. This gives us the positions of the crosses in the areaprime image. Once we draw the crosses we can easily obtain the sweep map image $\zeta(D)$ of $D$ by drawing the English partition whose removable corners are  the cells that contain the crosses. The  Blue squares in the area  cells  of $\zeta(D)$ are caused by the increasing diagonal
 pairs   $(3,5),(3,8),(5,8),(5,7),(1,7)$. In fact,  these   pairs of cars are precisely those producing the dinv of $PF$. We will show that column LLT's  generalize Dyck paths LLT's by constructing them from the  areaprime image of their parking functions. On the right in  the  above display we have 
the sweep map image $\zeta(D)$. The English partition above 
$\zeta(D)$ contains  $4$ crosses. If we remove the crosses created by  the  forced

\hfill
$
 \includegraphics[height=1.2 in]{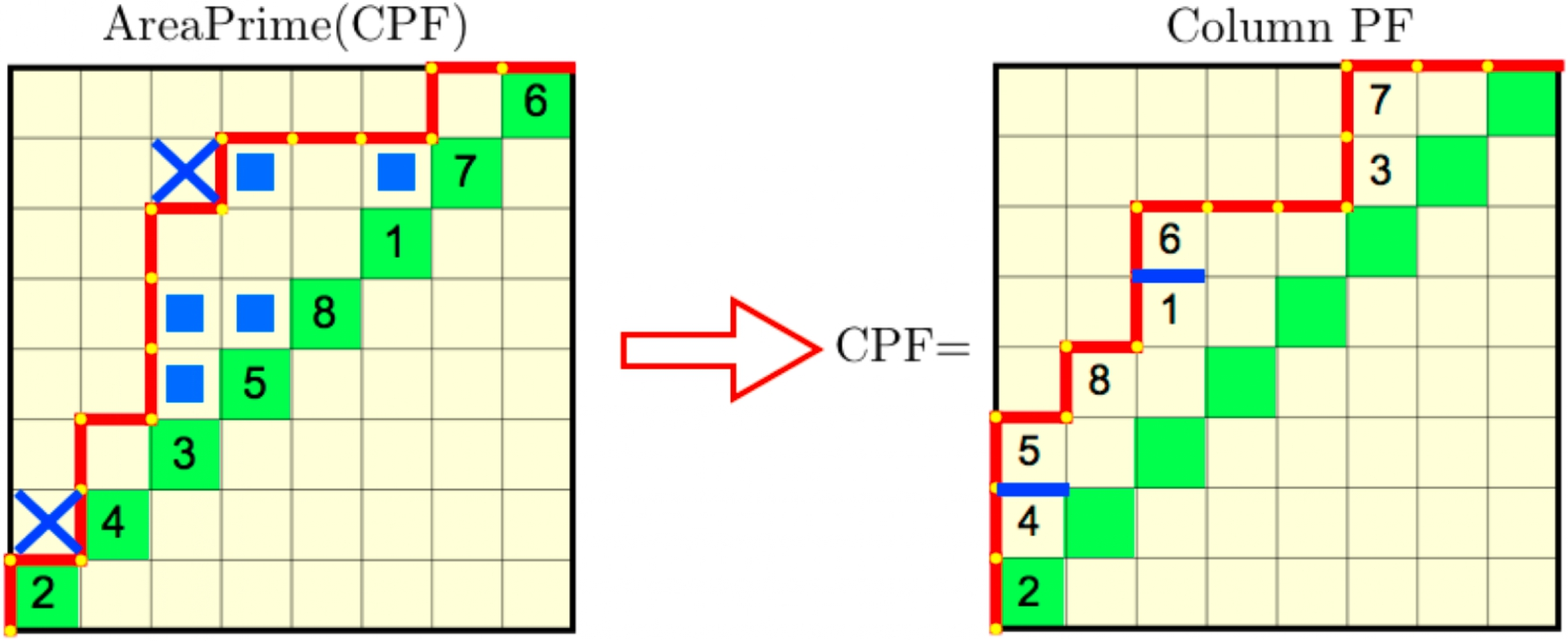}
$

\vskip -1.23 in
\hsize 3.4in
 
\noindent
 pairs $(4,5)$ and $(1,6)$ we are left with the  left portion of
this adjacent display. If we remove the $2$ remaining crosses and replace the North and East unit steps touching   each removed cross by a single  diagonal step the Dyck path $\zeta(D)$ becomes a Shroeder path. Since these diagonal steps can never occur on the   diagonal, the number

\hsize 6.5in
\noindent
of Shroeder paths thus obtained is one half of the Shroeder number.   In the next display we will transform this areaprime 
 with two removed crosses into what we will call a ``column parking function". 
   The figure on the right of the above display gives an intermediate step.
 \vfill
 \supereject
 
 \vskip .6in
 \hfill
$
 \includegraphics[height=1.2 in]{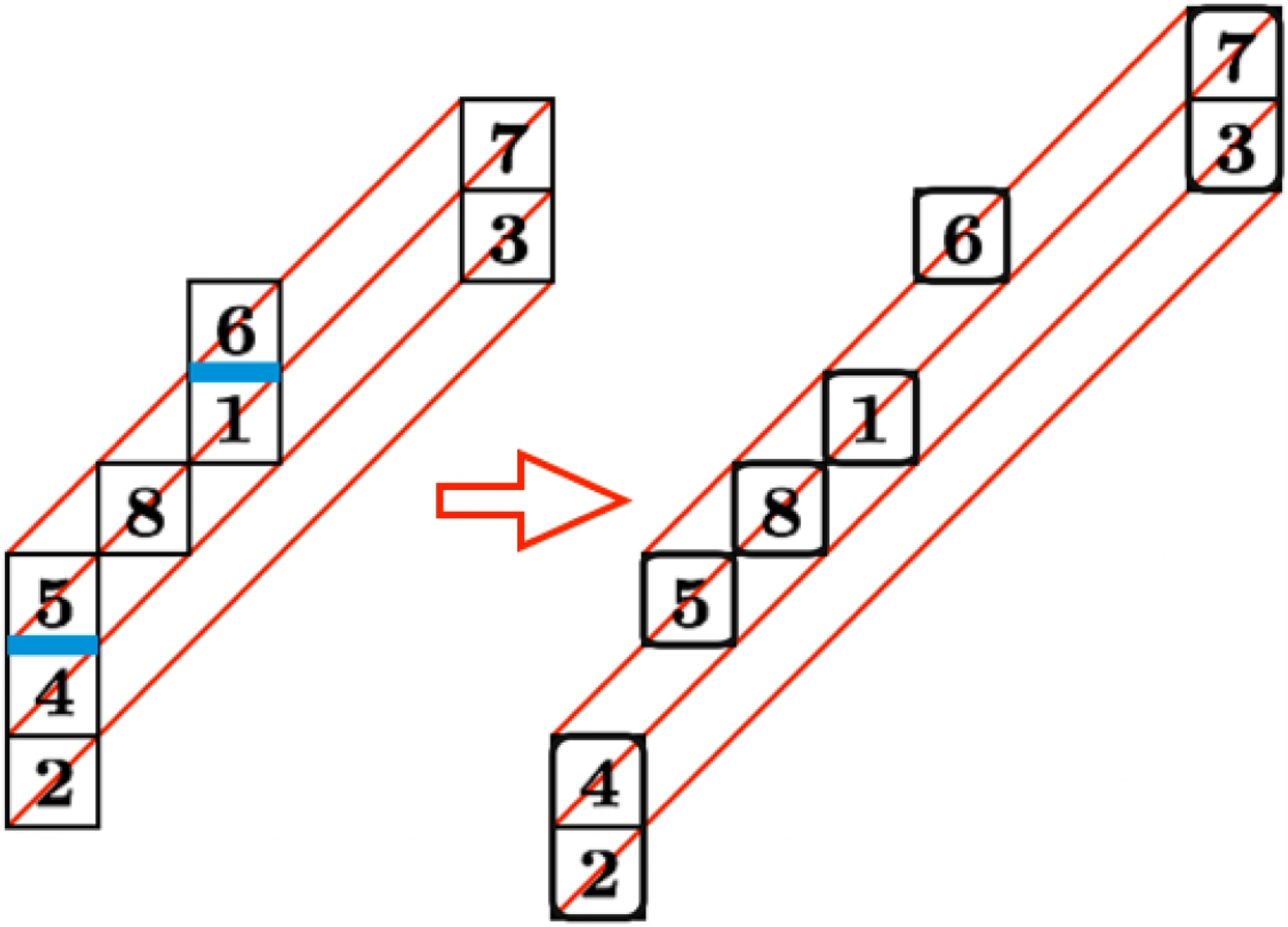}
$

\hsize 4.7in

 \vskip -1.4in  
\noindent
 To obtain it we simply insert in the original parking function, between cars $4,5$ and $1,6$, two separating blue dashes. In the left figure of
the adjacent display  we have simply  
reproduced only the columns with their blue dashes. In the right figure we have separated the cells containing $5$ and $6$  
from their columns. This  yields the right portion of the adjacent display. In summary, we have here identified an areaprime with  missing Blue crosses with a column parking function. In the display below we show how to construct the areaprime of a 

 \vskip -.04in
\hsize 6.5in
\noindent
column parking function. We have labeled  the  cars by letters but we assume that the 
cars they  represent are column increasing. 
 To construct its areaprime
we start by placing $\sig(CPF)$ in the diagonal as shown. Next we insert a blue dot at the center of any cell defined by a potential dinv. For instance
 the dot in the $3^{rd}$ row is the potential 
secondary dinv created when when  $b<c$. The dot in the $4^{th}$ row is due to the primary 

\hsize 4.2in

\vskip -.04in
\noindent
dinv caused by $c<d$. Every  one of the blue dots is caused by a potential dinv.
The next  step is to enclose all the blue dots by a Dyck path. 
Finally we add the blue crosses corresponding  to  the ``forced" pairs of cars. We see that not all the removable corners of the English partition above the path have crosses.  Confirming the fact that areaprimes of CPF    are none other than areaprimes of PF  with missing Blue crosses. 
  The LLT polynomial of a column LLT  may the be written in the form
\vskip -1.5in
\hsize 6.5in
\hfill
$
 \includegraphics[height=1.3 in]{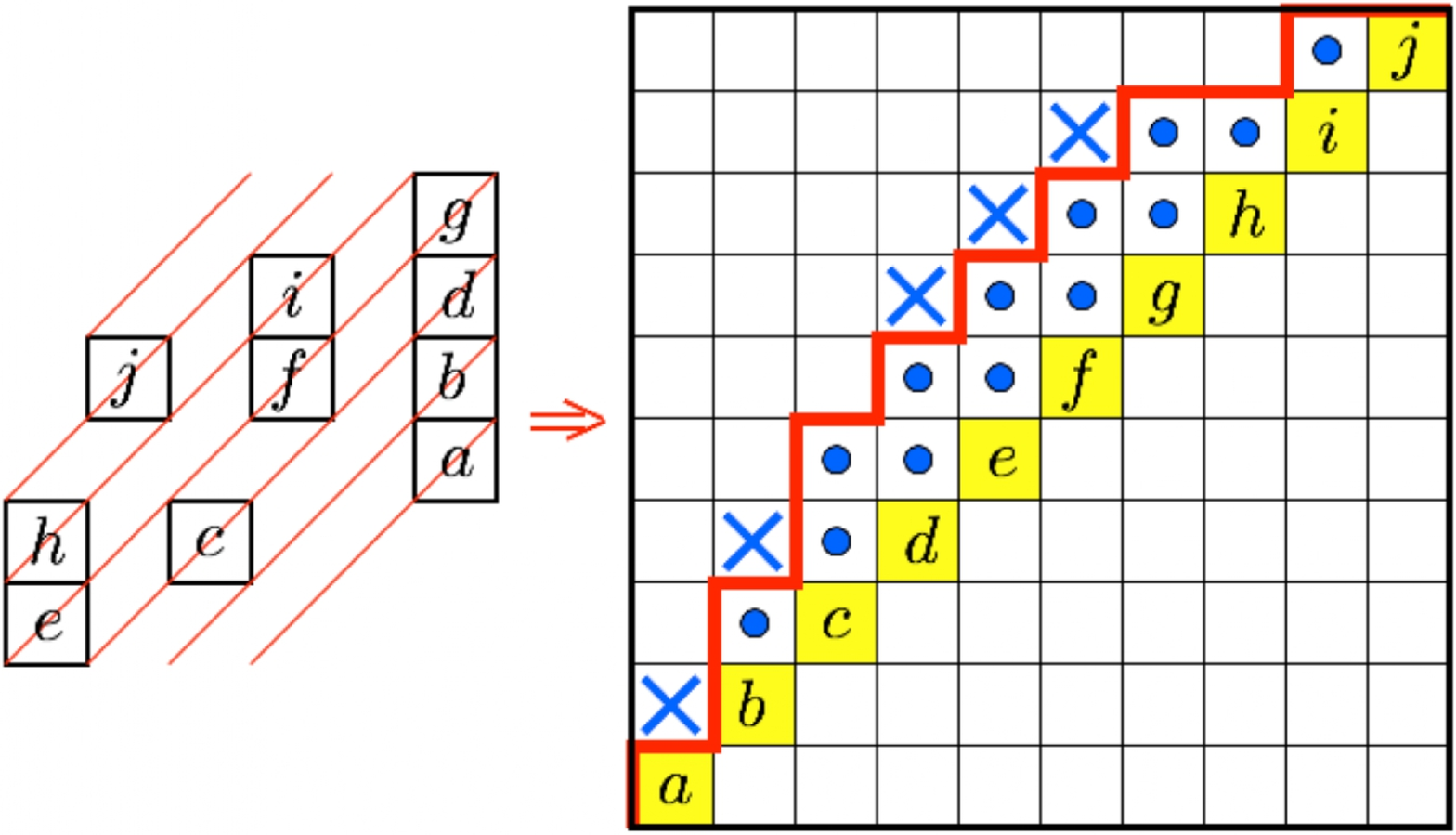}
$

\vskip .08in
\noindent

\vskip -.5in
$$
LLTC_{D,T}[X;q] \, =
\sum_{CPF\in {\cal CPF}_{T}}q^{dinv(CPF)}
F_{pides(CPF)}[X]
\eqno 3.15
$$
\vskip -.2in
\noindent
where $T\con forced(D)$, 
${\cal PF}_{T}$ is the  family of column parking functions  restricted to increase across  the pairs of cars in  $T$, and the composition $pides(CPF)$  as usual gives the descent set of
$\sig(CPF)^{-1}$. It also follows from our construction that the cardinality of these polynomials is given by the lower Schroeder number. 
 
\vskip -.06in
For later purposes it will be necessary to construct the same polynomial $LLTC_{D,T}[X;q]$ by starting from a general Dyck path  $Z=\zeta(D)$ in the $n \times n$ lattice square $R_n$ and a subset of pairs $T\con forced(D)$, by following the following rules.

\item{\bf 1)} {\ita Draw the path $Z$ where $Z[j]$ gives the number of coarea cells in the $j^{th}$ row of $R_n$.}  
\item{\bf 2)} { \ita A permutation $\sig\in S_n$  is called $Z,T$-compatible if and only if 
 for every pair $(r,s)\in T $ we have $\sig_{n+1-r}<\sig_{n+1-s} $.}
 
\item{\bf 3)} {\ita For a given $Z,T$-compatible $\sig$, a pair $1\le i<j\le n$
contributes  a unit of $dinv$   if and only if  $[i,j]$ is in the dinvset of $Z$ and $\sig_{n+1-j} >\sig_{n+1-i}$.}
\sas

\noindent
This given, we have the Schur expansion
\vskip -.2in
$$
LLTC_{D,T}[X;q]\, = \sum_{\sig\in \CF_{Z,T}}q^{dinv(\sig)}\, s_{pides(\sig)}[X].
\eqno 3.16
$$
\vskip -.2in
\noindent
where $\CF_{Z,T}$ is the family of all the $Z,T$-compatible $\sig\in S_n$, and  
`` $pides(\sig)$'' denotes the composition that gives the descent set of $\sig^{-1}$.

It turns out that this formula can be used only for moderately small $n\le 7$.
To confirm the validity of our $e$-positivity conjectural expansions, we will use the  Carlsson-Mellit super-fast manipulatorial way of computing the same symmetric polynomials. These new formulas have a complexity which is only linear in the number of steps of $Z$. We will present  them here in full detail since they    
are somewhat difficult to  extract out of the original paper [2].

\hfill
\supereject

These formulas are in terms of operators acting on the family $\Lambda[X;q,Y]$ of symmetric functions in the infinite alphabet $X$ with coefficients polynomials in $q$ and a variable alphabet $Y=\{y_1,y_2,y_3,\ldots\}$. The parameter  $k$ will denote the size of the $Y$-alphabet.

 The operators are $d_{+}$, $d_{-}$  and the bracket $[d_{-},d_{+}]$. They are all expressed in  terms of the operators $T_i$ with $T_0$ acting as identity and $T_i$ (for $i\ge 1$) acting on $F\in \Lambda[X;q,Y]$ precisely as follows
$$
T_i F\ses{ (q-1)y_iF\sps  (y_{i+1}-q y_i)\,s_iF \over y_{i+1}-y_i},
\eqno 3.17
$$
where $s_i$ denotes the transposition that interchanges $y_i$ with $y_{i+1}$.  
Using 3.17,  we set for $F\in \Lambda[X;q,y_1,\ldots,y_k]$

\hsize 4.8in

$$
d_{+}^k\, F\ses T_1T_2\cdots  T_k F[X+(q-1)y_{k+1}].
\eqno 3.18
$$
$$
\ess\ess\ess\ess\ess\ess
d_{-}^k\, F\,=\, -  F[X-(q-1)y_{k}]\sum_{i\ge 0}(-1/y_k)^ie_i\Big|_{y_k^{-1}}.
\eqno 3.19
$$
$$
[d_{-},d_{+}]^k \,F\, =\, {d_{-}^{k+1}d_{+}^kF\sms d_{+}^{k-1}d_{-}^{k}F\over q-1}.
\eqno 3.20
$$
  
\hsize 6.5in

\vskip -1.36in
\hfill
$ \includegraphics[height=1.3 in]{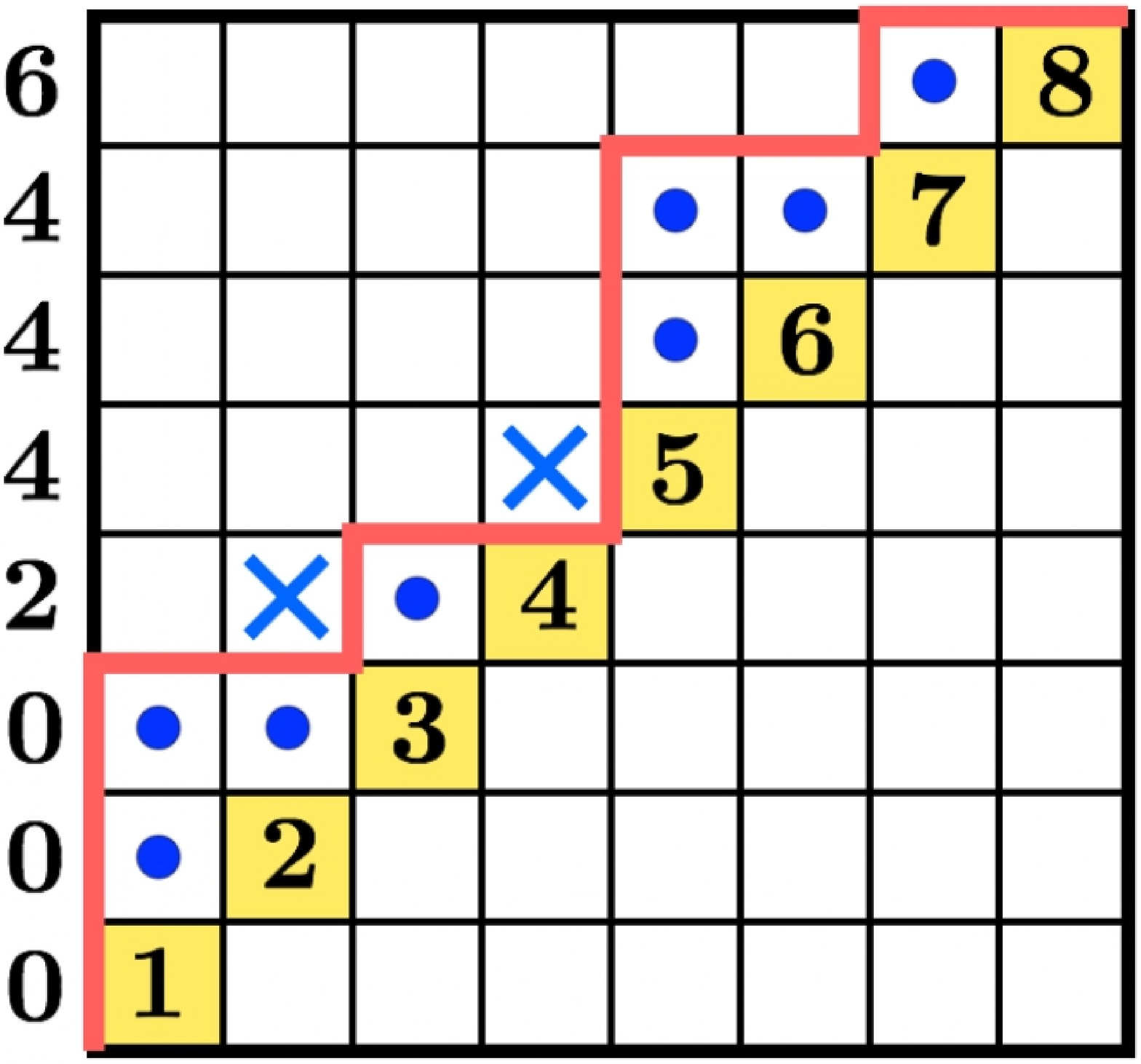}  $

These operators are used in a very simple manner to obtain Dyck path LLT's, column  LLT's and unicellular LLT's. We need only carry out the details in a special case. In the display above we have the areaprime image  of a typical  column LLT. Here $Z=[0,0,0,2,4,4,4,6]$. Notice that the English  partition above $Z$ has three removable corners $[[2,4],[4,5],[6,8]]$ but only the first two are marked. So in this case $T=[[2,4],[4,5]]$. The permutation we placed in the diagonal corresponds to the maximal column parking function, the one whose $dinv$ is the number of blue dots under $Z$. Using a ``$0$'' for each North step and a ``$1$'' for each East step, the sequence of steps of $Z$ can be represented by the following word
$$
 \includegraphics[height=.3 in]{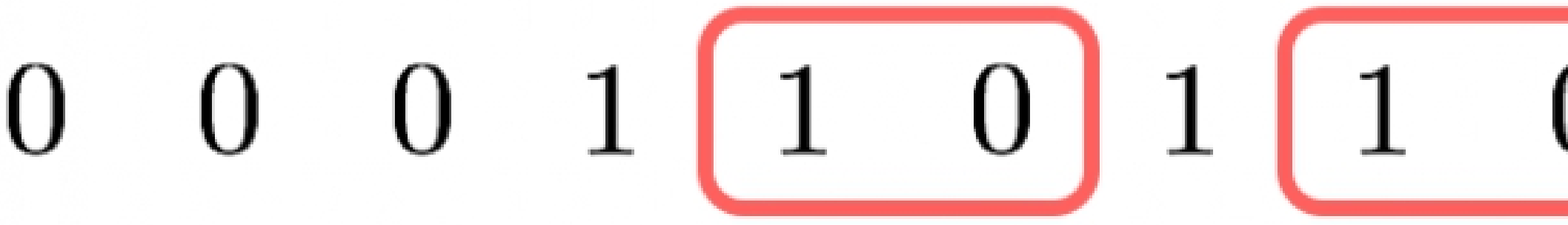}   \ ^\cdot
\eqno 3.21
$$
We can immediately identify the three removable corners by simply locating the East steps followed by a North step. We purposely framed the marked ones. The final word we will use to  guide the construction of the polynomial $LLTC_{D,T}[X;q]$ where $Z=\zeta(D)$ is obtained by replacing the framed pairs by a ``$2$'' obtaining the compressed word 
$$
 \includegraphics[height=.17 in]{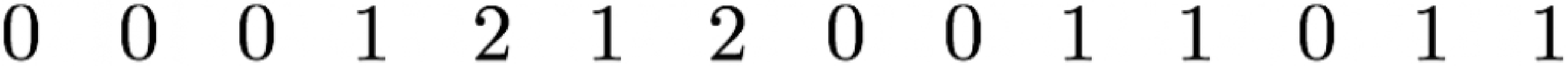} \ \cdot
\eqno 3.22
$$
Starting with the symmetric function $F=1$, and proceeding from right to left,
we apply a ``$d_{+}$'' for each $1$, a ``$d_{-}$'' for each $0$, and a ``$[d_{-},d_{+}]$'' for each $2$, according to the following sequence of commands 
based on $W$ being the reverse of the word in 3.22:

\itemitem{1)} {\ita set $k=0$; set $out=1$;}
\itemitem{2)} {\ita for i from 1 to 14 do}
\itemitem{3)} {\ita if W[i]=1 then $out=d_{+}^kout;\ k=k+1; $}
\itemitem{4)} {\ita else if W[i]=2 then $out=[d_{-},d_{+}]^k out$;}
\itemitem{5)} {\ita else $out=d_{-}^k out;\ k=k-1; $}
\itemitem{6)}{\ita end if;$\,$ end do; }
  
The general case is easy to derive from this example. What is remarkable about this algorithm,  is not only that it is of linear complexity but that it can  be used in all  three  cases: Dyck path LLT's ($T=forced(Z)$, all corners marked), column LLT's ($T\con forced(Z)$, some corners marked), unicellular LLT's ($T=\phi$, no corners marked).

We will see later how this Carlsson-Mellit way of obtaining the column LLT polynomials  can be used to check, for relatively large $n$, the recursive way of constructing our conjectured $e$-expansions.   

\vfill\supereject

\hsize 6.5 in
\noindent

\centerline{\ita b) The extension to column  LLT's of  Conjecture 3.1.} 
\sas

\noindent{\bol Conjecture 3.2}

{ Given a Dyck path $D$ in the $n\times n$ lattice square $R_n$ and a  subset 
$T\con forced(D)$, the following algorithm 
constructs the $e$-basis expansion of the polynomial $LLTC_{D,T}[X;1+q]$.  Draw the parking function $CPF$ with $\sig(CPF)=n\cdots 321$ and determine the set of pairs of cars $(a,b)$ with $a<b$ producing a dinv. Call this ``$dinvset(D)$''. 
Then our final product can be written  in the form

\vskip -.15in
$$
LLTC_{D,T}[X;1+q] \, =
\sum_{S\con dinvset(D)}q^{|S|}
e_{\mu(S,T)},
\eqno 3.16
$$

\vskip -.15in
\noindent
where the polynomial
$$
\sum_{S\con dinvset(D)}q^{|S|}
e_{\mu(S,T)}
$$
is obtained by repetitions of the $4$ steps stated in Conjecture 3.1 except that in the second step the ``$forced (D)$" set of pairs  must be replaced by $T$.

\centerline{\ita c) The conjecture that extends to all column  LLT's the recursion in Theorem 3.2}
\sas

Our point of departure, as before,
 is a Dyck path $D$ in $R_n$. However, here it will be
convenient to work entirely with the areaprime image. Thus we let $Z=\zeta(D)$ and fill, from bottom to top, all the diagonal   

\hfill
$
 \includegraphics[height=1.5 in]{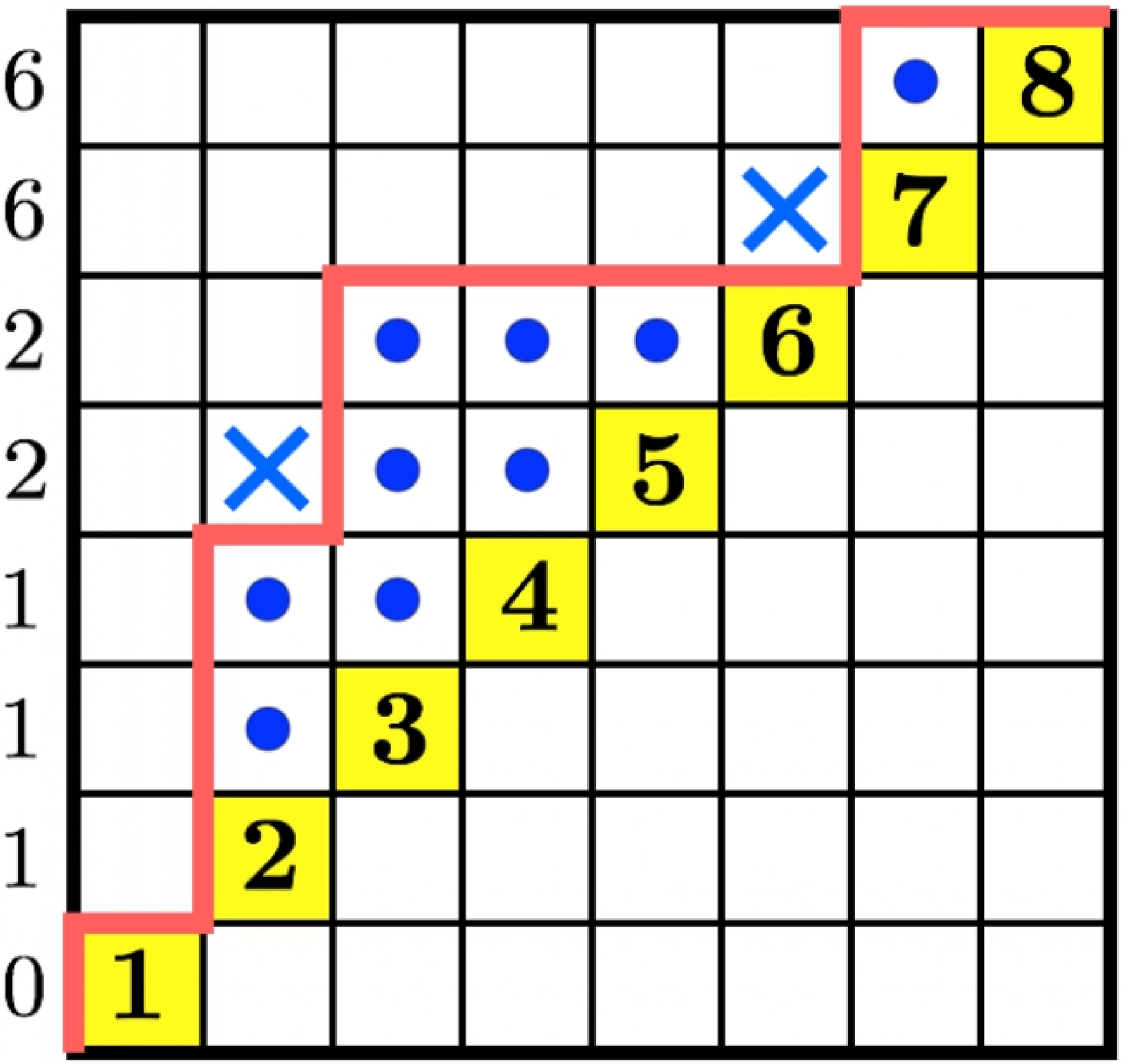}
$

\vskip -1.55in
\hsize 4.8in
\noindent
cells with ``cars'' $1,2,\ldots,n$. Again, it will suffice to carry out our algorithm in a special case. 
Our choice is $Z=[0,1,1,1,2,2,6,6]$, we also choose 
\hbox{$T=[[2,5],[6,7]]$.} These two pairs are represented by the blue crosses in our display. The blue dots below $Z$ represent, in the same  manner, the pairs in $dinvset(Z)$.  The first step in the algorithm is to determine the {\ita possible}  downsets of $n$. In our case, $n=8$. This requires two important properties. A subset $S=\{ i_1<i_2<\cdots <i_k=n\} \subseteq \{1,\ldots,n\}$ is ``{\ita possible}'' if and only if every pair of elements $(i_j,i_{j+1})$ in $S$ is in the union of pairs $T \cup dinvset(Z)$; moreover,  whenever  $(i,j)\in T$ and  $j\in S$ 

\hsize 6.5in
\noindent 
 then $S$ must also contain $i$. These two properties are easy consequences of the manner our downsets are constructed. Our implementation of this algorithm yields only $12$ possible downsets. The next step of our algorithm is to compute the weight of each {\ita possible} subset $S$ of cars. This consists of   a polynomial $weight_S(q)$. The computation of this weight is the most delicate step of our algorithm.
 We  illustrate our procedures by displaying the output in the two cases $S_1=\{2, 3,5, 6,7,8\}$ and $S_2=\{ 3,4,6,7,8\}$.    

\vskip .02in
\noindent
$
\ess\ess\ess \includegraphics[height=1in]{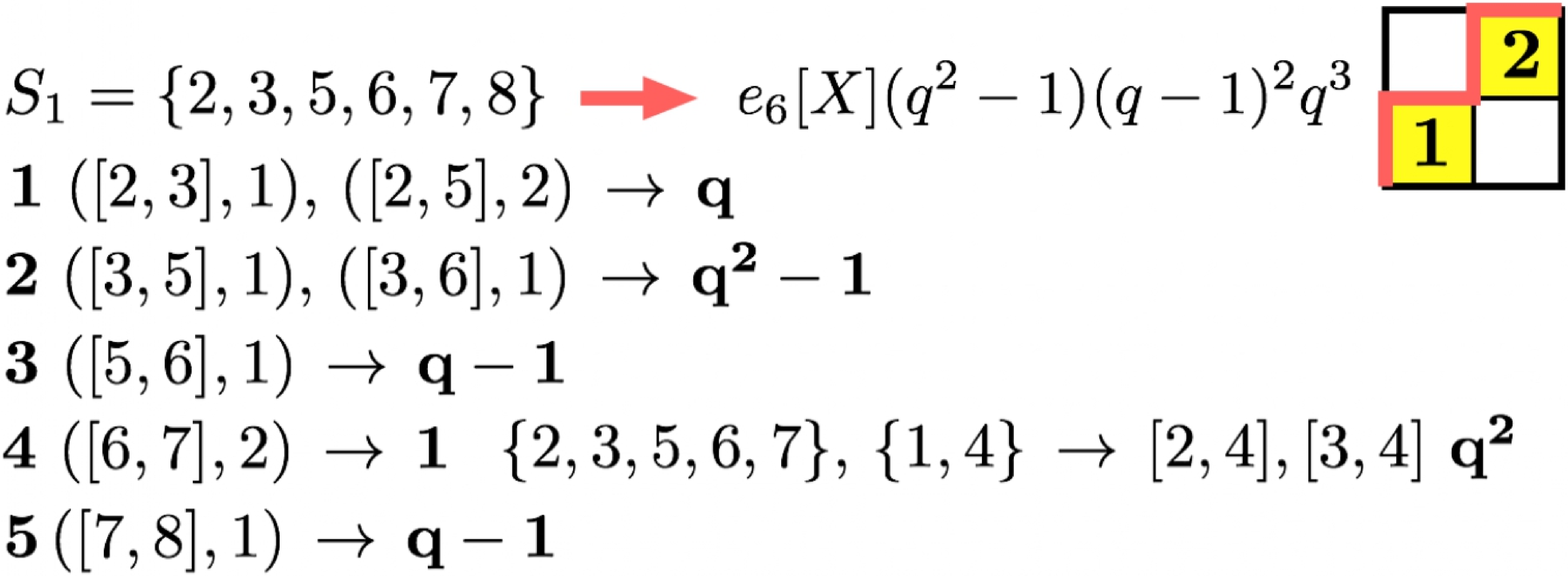} 
$

\vskip -1.03in
\hfill
$
\includegraphics[height=1in]{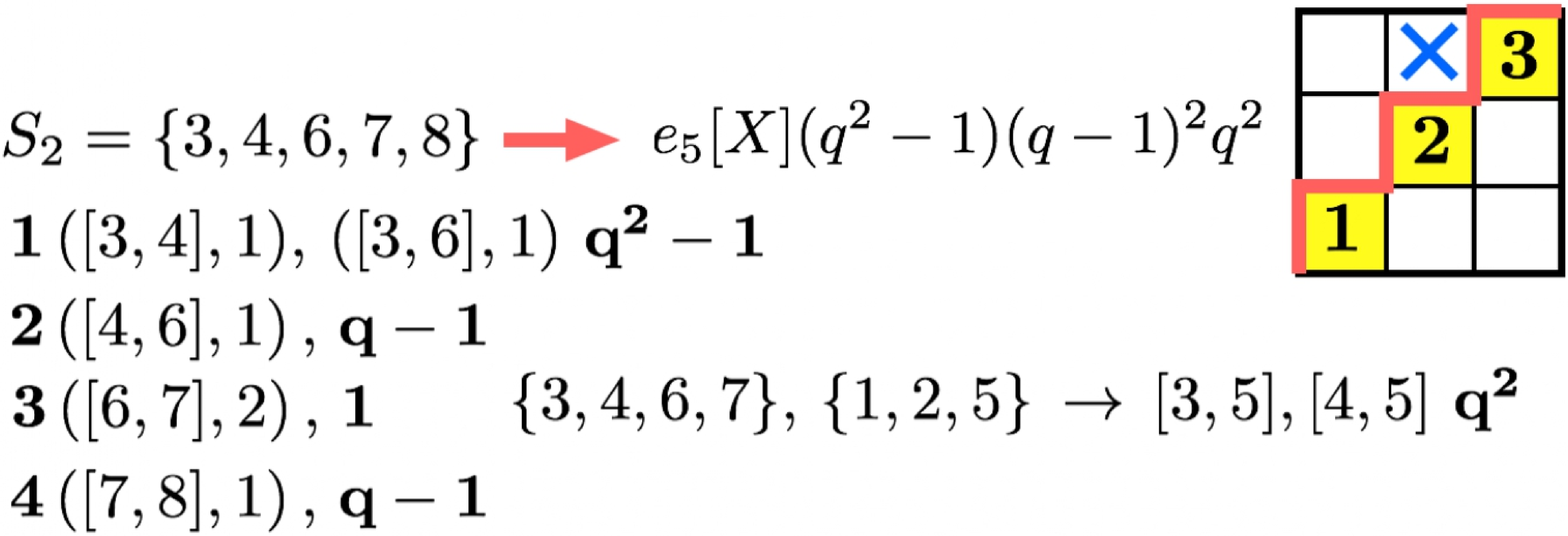}\ess 
$

\vskip -.02in 
\noindent
Our procedures yield  the conjectured $e$-expansion of the polynomial 
$LLTC_{D,T}[X;q]$ in terms of $Z=\zeta(D)$ and the pairs in $T$. The basic principle is to  use the split  
$q=q-1+1$ and take $\bf q-1$ or $\bf 1$ as the case may be. We take both when the result does not affect  $downset(n)$, this yields a contribution of $\bf q$ to its weight.

For instance, $S_1-\{8\}=\{2,3,5,6,7\}$ that leaves $\{1,4\}$, Since the additional edges $[2,4]$ and $[3,4]$ do not affect $downset(8)$ we can add $\bf q^2$ to its weight. The analogous argument holds for $S_2$, with the additional edges $[3,5]$ and $[4,5]$, as indicated in the adjacent display. To obtain the total weight of $downset(8)$ for $S_1$, each of its elements are processed one at the time, as indicated in the left display. In the column of $2$ there is only a blue dot in the row of $3$ that accounts for the $1$ in $([2,3],1)$. In the same column we see a blue cross in the row of $5$ (see above display). This accounts for the $2$ in $([2,5],2)$. Now since blue crosses do not contribute the weight, the weight of the column of $2$ reduces to $\bf q$. The reason for this is that we may or may not add the edge $[2,3]$. When we process the column of $3$ whether we add at least one of  the edges  $[3,5]$ or 

\vfill
\supereject
 
 \noindent
$[3,6]$ by connecting $3$ to $5$ or $6$ we will guarantee the addition of $3$ to $downset(8)$. This accounts for the contribution of a $\bf q^2-1$ of the column of $3$. The reader should  have no difficulty interpreting the outputs of the two cases in the previous display.
 
 Once all the possible {\ita downsets} are processed, we  obtain
 \vskip -.15in
$$
areaprime_{Z,T}[X,q]\,= \bu \bu\sum_{S\,\in\, possible_{Z,T}}\bu\bu e_{|S|}[X]\, weight_S(q)\,
areaprime_{Z(S),T(S)}[X,q],
\eqno 2.17
$$ 
\vskip -.1in
\noindent
where the Dyck path $Z(S)$ and the residual $T(S)=forced \big(Z(S)\big)$ are both obtained by simply deleting each element of $S$ in the diagonal of $R_n$ along with every cell of $R_n$ that is in the same row and column as that element. Notice that the validity of our conjectured $e$-expansion can now be checked for relatively large examples by computing all symmetric polynomials $areaprime_{Z,T}[X,q],$ and $areaprime_{Z(S),T(S)}[X,q]$ by the fast algorithm of
Carlsson-Mellit.

\vskip -.1in
\hfill
$
 \includegraphics[height=1in]{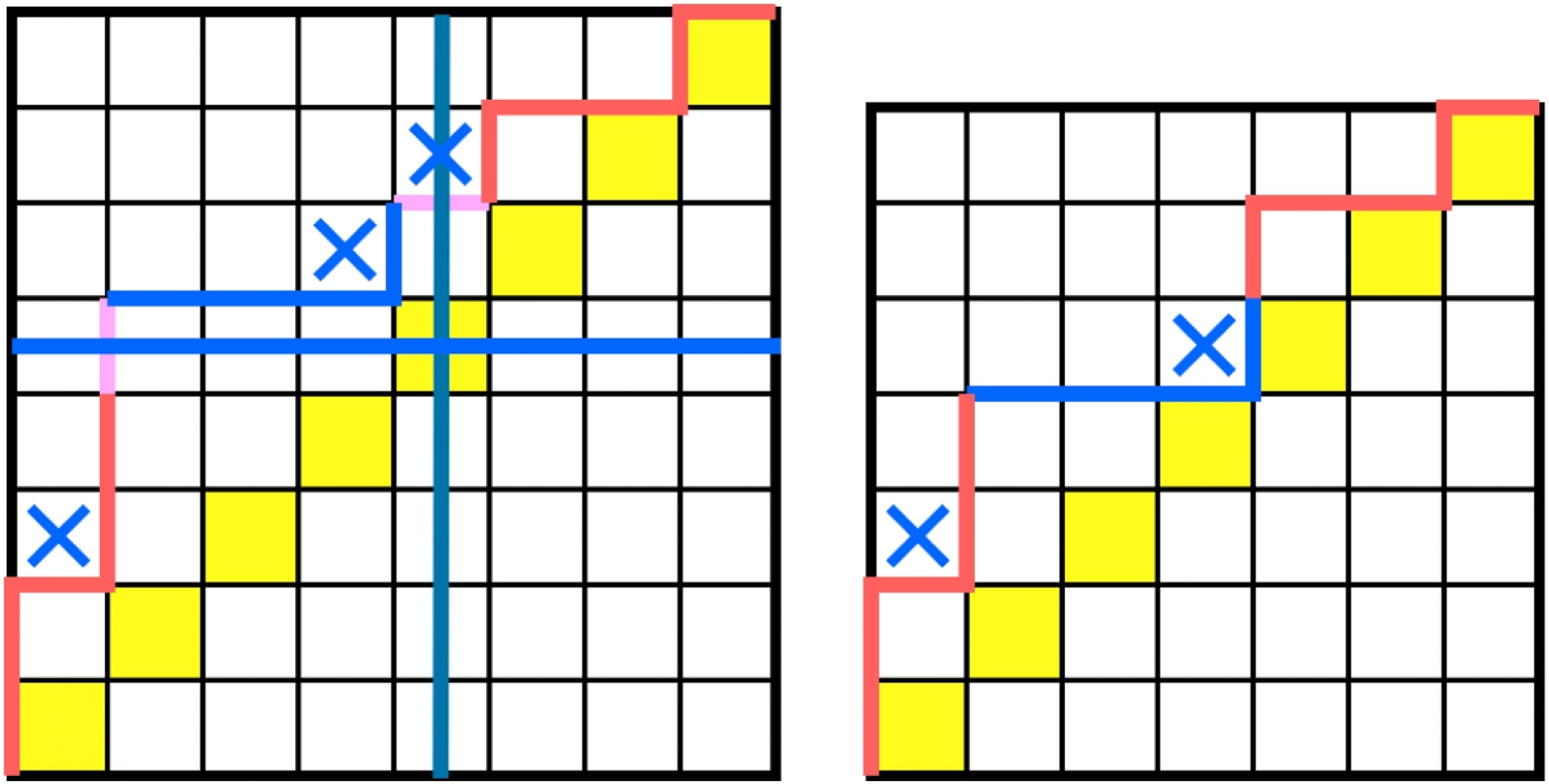}
$

\vskip -.94in
\hsize 4.4in
It may be good to see that the deletions of rows and columns invariably yield that all the resulting $Z(S)$ 
are Dyck paths.  By induction the following argument should be sufficient.
Notice first that since the row and column we are deleting meet at a diagonal cell, except in the case this cell is the first or the last, the Dyck path $Z$ is broken up into 

\hsize 6.5in
\noindent
three pieces. Leaving aside that limit  case, in the left of the above display, we have depicted the first and the last of these three pieces in red, the middle piece in blue and the two deleted steps in pink.
Now the new path $nZ$ consists of all the steps of $Z$  except the two  deleted steps. To construct the resulting path $nZ$, we start by the steps of $Z$ up to the deleted North step. That is the first red piece of $Z$ in our display. Now $nZ$
will continue by the steps of the middle piece, in blue in our display. To compensate for the deleted North step, the middle piece of $Z$ must move vertically down one unit.  Next comes the deleted East step. To compensate for the deleted steps $nZ$ must follow
the steps of $Z$ shifted diagonally by one cell.
To prove that $nZ$ is also a Dyck path, we need only show that it remains weakly above the diagonal. Since first piece of $Z$ and $nZ$ are identical and the last piece of $Z$ moves diagonally, neither can cross the diagonal. 
However, the middle part of $nZ$ cannot cross the 
diagonal either. In fact,  the middle part  of $Z$ must all be strictly above the diagonal, due to the center of deletion being a diagonal cell. 

The marginal cases we left aside can be dealt with in the same  way since in the first case the first piece is missing and in the second case  the last piece is missing.  The remaining pieces can be dealt with exactly as in the above argument.
\sas

It must be mentioned  the  finding of Kreweras in [15] that hit the tip of an iceberg is the surprising relation between $\nabla e_n$ and the Kreweras polynomials defined in [15] as the family satisfying the recursion and base case 
$$
P_{n+1}(q)\ses \sum_{i=0}^n{n \choose i}[i+1]_q P_i(q)P_{n-i}(q),
\bigsp \ess\ess\hbox{$P_0(q)=1$}.
$$ 
The precise relation between the Kreweras polynomials and $\nabla e_n$, that follows from our conjectured  $e$-positivity phenomenon 
of Dyck path LLT's, is
$$
P_{n}(1+q)\ses \del_{p_1}^n\big(\nabla e_n[X;1,1+q]\big).
$$

\sas
It is clear that all these $e$-positivities  create  a variety of new problems in Algebraic Combinatorics. We might even say that the Dyck paths LLT conjecture by itself creates a major upheaval in this field.

\vfill\supereject

\centerline{\bol Bibliography}
\sas

\item {[1]} P. Alexandersson, {\ita LLT Polynomials, Elementary Symmmetric Functions and Melting Lollypops}, 

$\ssp  $ arXiv:1903.03998 (2019).
\sa

\item {[2]} E. Carlsson and A. Mellit,
{\ita  A Proof of the Shuffle Conjecture}, arXiv:1508.06239 (2015).
\sa

\item {[3]}
E. Egge, N. Loehr and G. Warrington, 
{\ita From quasisymmetric expansion to Schur expansion via 
a modified inverse Kostka matrix}, European J. Combin., 
{\bf 31 (8)} (2010), 2014--2027.
\sa

\item  {[4]} A. M. Garsia, G. Xin and M. Zabrocki, {\ita  Hall-Littlewood operators in the Theory of Parking Functions and Diagonal Harmonics},
International  Mathematical Research Notices, {\bf 2012 (6)} (2012), 1264--1299.
\sa

\item {[5]} A. M. Garsia and M. Haiman,
{\ita A Remarkable  $q,t$-Catalan sequence and  $q$-Lagrange Inversion},
 J. Algebraic Combin., {\bf 5 (3)} (1996), 191--244.
\sa

\item {[6]} A. M. Garsia and J. Haglund,
{\ita A proof of the $q,t$-Catalan positivity conjecture},
Discrete Mathematics, {\bf 256 (2)} (2002), 677--717.
\sa

\item {[7]} A. M. Garsia, J. Haglund, J. B. Remmel and M. Yoo, {\ita A proof of the Delta Conjecture at q=0}, arXiv:1710.07078 (2017).
\sa

\item {[8]} A. M. Garsia and J. B. Remmel, {\ita 
A note on passing from a quasi-symmetric expansion to a Schur expansion of a symmetric function}, arXiv:1802.09686 (2018).
\sa

\item {[9]}I. Gessel, {\ita Multipartite $P$-partitions and inner products of skew Schur functions}, Contemp. Math., {\bf 34} (1984), 289--301.
\sa

\item  {[10]} I. Gessel and D. Wang, {\ita  Depth-first search as a combinatorial correspondence}, J. Combin. Theory Series. A, {\bf 26 (3)} (1979), 308--313.
\sa

\item {[11]} J. Haglund, M. Haiman, N. Loehr, J. B. Remmel and A. Ulyanov,
{\ita A combinatorial formula for the character of the diagonal coinvariants},
{ Duke Math. J.}, {\bf 126} (2005), 195--232.
\sa

\item  {[12]}  J. Haglund, {\ita The q,t-Catalan numbers and the space of diagonal harmonics}, volume 41 of
University Lecture Series. American Mathematical Society, Providence, RI (2008).
\sa

\item  {[13]} J. Haglund,  J. Morse and M. Zabrocki,
{\ita
A Compositional Shuffle Conjecture Specifying Touch Points of the Dyck Path}.
Canad. J. Math., {\bf  64 (4)} (2012), 822--844.
\sa

\item  {[14]}  J. Haglund, J. B. Remmel and A. Wilson, {\ita  The Delta Conjecture}, arXiv: 1509.07058 (2017).
\sa

\item  {[15]} G. Kreweras, {\ita Une famille de polynomes ayant plusieurs proprietes enumeratives}, Periodica
Math. Hungar., {\bf 11} (1980), 309--320.
\sa

\item {[16]}  I. G. Macdonald.
{\ita Symmetric functions and Hall polynomials},
 2nd Ed. Reprint of the 2008 paperback edition,
 Oxford University Press, New York (2015).
 \sa
 
\item {[17]}  J. Novak, {\ita  Three lectures on free probability}, arXiv:1205.2097 (2012).
\sa

\item {[18]}  M. Zabrocki,  {\ita  A module for the Delta conjecture},  arXiv:1902.08966 (2019).

 \end